# Output-Feedback Control of Viscous Liquid-Tank System and its Numerical Approximation

**Iasson Karafyllis[*], Filippos Vokos[*] and Miroslav Krstic[**]**

[*]Dept. of Mathematics, National Technical University of Athens, Zografou Campus, 15780, Athens, Greece,
emails: iasonkar@central.ntua.gr, fivo33@hotmail.com

[**]Dept. of Mechanical and Aerospace Eng., University of California, San Diego, La Jolla, CA 92093-0411, U.S.A., email: krstic@ucsd.edu

**Abstract**

We solve the output-feedback stabilization problem for a tank with a liquid modeled by the viscous Saint-Venant PDE system. The control input is the acceleration of the tank and a Control Lyapunov Functional methodology is used. The measurements are the tank position and the liquid level at the tank walls. The control scheme is a combination of a state feedback law with functional observers for the tank velocity and the liquid momentum. Four different types of output feedback stabilizers are proposed. A full-order observer and a reduced-order observer are used in order to estimate the tank velocity while the unmeasured liquid momentum is either estimated by using an appropriate scalar filter or is ignored. The reduced order observer differs from the full order observer because it omits the estimation of the measured tank position. Exponential convergence of the closed-loop system to the desired equilibrium point is achieved in each case. An algorithm is provided that guarantees that a robotic arm can move a glass of water to a pre-specified position no matter how full the glass is, without spilling water out of the glass, without residual end point sloshing and without measuring the water momentum and the glass velocity. Finally, the efficiency of the proposed output feedback laws is validated by numerical examples, obtained by using a simple finite-difference numerical scheme. The properties of the proposed, explicit, finite-difference scheme are determined.

**Keywords:** Saint-Venant model, shallow water equations, output feedback, Control Lyapunov Functional, PDEs.

## 1. Introduction

Output feedback stabilization of Partial Differential Equations (PDEs) is an active research area. Important results in the literature include the construction of observer-based output feedback controllers for parabolic or hyperbolic PDEs (see for instance [15,27,28,29,38,44,45]). Cases of systems with boundary measurements subject to time-delays are also considered in the literature (see for instance [23,24,37]).

The major control scheme in the linear case is the so-called "Separation Principle" (see [39]) for the finite-dimensional case): a state feedback controller is designed separately from an observer which estimates the whole state of the PDE system. However, this scheme is difficult to be applied in the nonlinear case because there are no systematic procedures for the design of observers which estimate the infinite-dimensional state of the system. For complex systems it may be reasonable to set a less ambitious goal for output feedback stabilization. Luenberger in [32,33] introduced the concept of a functional observer: a functional observer is an auxiliary system that is driven by the input and available system output and estimates a function of the state. Research on functional observers continues (see [42] and references therein). As already remarked by Luenberger in [33], it is possible to utilize a functional observer in order to directly estimate the nominal state feedback (see also [25]).

In this paper, we apply the functional-observer methodology to a well-known nonlinear PDE problem: the stabilization of the motion of a tank containing a liquid described by the viscous Saint-Venant PDE system. The Saint-Venant system or shallow water model has been used in many applications (river flows, tidal waves, oceans) and since its first derivation from first principles by Adhémar Jean Claude Barré de Saint-Venant in 1871 (see [1]), the model has been extended to take into account many types of forces acting on the body of the liquid other than gravity (e.g., viscous stresses, surface tension, friction forces; see [5,6,16,30,34,41]). Feedback stabilization problems involving various variations of the Saint-Venant model have attracted the attention of many researchers during the last decades (see [2,3,4,7,8,9,10,12,13,21,31,35,36,38]). Most works study the inviscid Saint-Venant model (i.e., the model that ignores viscous stresses and surface tension and takes into account gravity and friction forces), which provides a system of first-order hyperbolic Partial Differential Equations (PDEs). Many papers study stabilization problems for the linearization of this system around an equilibrium point, employing either the backstepping methodology (see [12,13]) or the Control Lyapunov Functional (CLF) methodology (providing local stabilization results for the original nonlinear system in many cases; see [2,3,4,8,9,10,31]). In [21] a state feedback stabilizer for the viscous, nonlinear, Saint-Venant system has been proposed. The feedback law proposed in [21] requires measurement of the tank position and velocity, the liquid level at the tank walls and the total liquid momentum. The latter quantity is very difficult to be measured directly.

In this paper we utilize the state feedback law proposed in [21] in order to solve the output feedback stabilization problem of the motion of a tank containing a viscous, incompressible, Newtonian liquid. We assume that we measure the tank position and the liquid level at the tank walls and we construct functional observers for the tank velocity and the liquid momentum. The problem is studied by means of the CLF methodology initially presented in [26] for nonlinear parabolic PDEs and studied also in [8,18,19,20,21]. We construct two different observers for the tank velocity: a full-order observer and a reduced-order observer. We also construct two different estimators for the liquid momentum (strictly speaking they are not functional observers of the liquid momentum). Therefore, we construct four different types of output feedback stabilizers. Exponential convergence of the closed-loop system to the desired equilibrium point is achieved in each case. As in [21], the constructed CLFs are also size functionals (see [40] for the notion of the size function) as they provide positive upper and lower bounds of the liquid level. As far as we know, this is the first paper in the literature that achieves output feedback stabilization of the nonlinear viscous Saint-Venant system.

The paper is structured as follows. In Section 2 we present the output feedback control problem of the motion of a tank. In Section 3 we present the main ideas for the construction of dynamic output feedback laws. Section 4 contains the statements of the main results. In Section 5 we present an algorithm which guarantees that a robotic arm will move a glass of water to a pre-specified position. We show that *no matter how full the glass is, the robot can transfer the glass without spilling water out of the glass, without residual end point sloshing and without measuring the water momentum* *and the glass velocity*. Section 6 presents illustrative numerical examples which show the efficiency of the proposed dynamic output feedback laws. While there are many numerical schemes devoted to the numerical approximation of the inviscid Saint-Venant model, there are only few papers in the literature that study the numerical approximation of the viscous Saint-Venant system with no control (see [11,17,22]) and there is no paper that studies the numerical approximation of the viscous Saint-Venant model with control. Here, we used a simple finite-difference scheme and in Section 6 we state the properties and the accuracy of the scheme. The proofs of all results are provided in Section 7. Finally, the concluding remarks are given in Sections 8.

**Notation.** Throughout this paper, we adopt the following notation.

* $\Re_+ = [0,+\infty)$ denotes the set of non-negative real numbers.

* Let $S \subseteq \Re^n$ be an open set and let $A \subseteq \Re^n$ be a set that satisfies $S \subseteq A \subseteq cl(S)$. By $C^0(A;\Omega)$, we denote the class of continuous functions on $A$, which take values in $\Omega \subseteq \Re^m$. By $C^k(A;\Omega)$, where $k \geq 1$ is an integer, we denote the class of functions on $A \subseteq \Re^n$, which takes values in $\Omega \subseteq \Re^m$ and has continuous derivatives of order $k$. In other words, the functions of class $C^k(A;\Omega)$ are the functions which have continuous derivatives of order $k$ in $S = \text{int}(A)$ that can be continued continuously to all points in $\partial S \cap A$. When $\Omega = \Re$ then we write $C^0(A)$ or $C^k(A)$. When $I \subseteq \Re$ is an interval and $G \in C^1(I)$ is a function of a single variable, $G'(h)$ denotes the derivative with respect to $h \in I$.

* Let $I \subseteq \Re$ be an interval, let $a < b$ be given constants and let $u : I \times [a,b] \to \Re$ be a given function. We use the notation $u[t]$ to denote the profile at certain $t \in I$, i.e., $(u[t])(x) = u(t,x)$ for all $x \in [a,b]$. When $u(t,x)$ is (twice) differentiable with respect to $x \in [a,b]$, we use the notation $u_x(t,x)$ ($u_{xx}(t,x)$) for the (second) derivative of $u$ with respect to $x \in [a,b]$, i.e., $u_x(t,x) = \frac{\partial u}{\partial x}(t,x)$ ($u_{xx}(t,x) = \frac{\partial^2 u}{\partial x^2}(t,x)$). When $u(t,x)$ is differentiable with respect to $t$, we use the notation $u_t(t,x)$ for the derivative of $u$ with respect to $t$, i.e., $u_t(t,x) = \frac{\partial u}{\partial t}(t,x)$.

* Given a set $U \subseteq \Re^n$, $\chi_U$ denotes the characteristic function of $U$, i.e. the function defined by $\chi_U(x) := 1$ for all $x \in U$ and $\chi_U(x) := 0$ for all $x \notin U$. The sign function $\text{sgn} : \Re \to \Re$ is the function defined by the relations $\text{sgn}(x) = 1$ for $x > 0$, $\text{sgn}(0) = 0$ and $\text{sgn}(x) = -1$ for $x < 0$.

* Let $a < b$ be given constants. For $p \in [1,+\infty)$, $L^p(a,b)$ is the set of equivalence classes of Lebesgue measurable functions $u : (a,b) \to \Re$ with $\|u\|_p := \left( \int_a^b |u(x)|^p \, dx \right)^{1/p} < +\infty$. $L^\infty(a,b)$ is the set of equivalence classes of Lebesgue measurable functions $u : (a,b) \to \Re$ with $\|u\|_\infty := \text{ess}\sup_{x \in (a,b)} (|u(x)|) < +\infty$. For an integer $k \geq 1$, $H^k(a,b)$ denotes the Sobolev space of functions in $L^2(a,b)$ with all its weak derivatives up to order $k \geq 1$ in $L^2(a,b)$.

## 2. Description of the Problem

We study a 1-D model for the motion of a tank that contains a viscous, Newtonian, incompressible liquid. We assume that a force, that can be manipulated, acts on the tank. Assuming that the liquid pressure is hydrostatic, the liquid is modeled by the 1-D viscous Saint-Venant equations and the tank obeys Newton's second law.

The control objective is to drive asymptotically the tank to a specified position without liquid spilling out and having both the tank and the liquid within the tank at rest. Figure 1 shows a picture of the problem.

Let the position of the left wall of the tank at time $t \geq 0$ be $a(t)$ and let the length of the tank be $L > 0$ (a constant). The following equations describe the motion of the liquid within the tank:

$$H_t + (HV)_z = 0, \text{ for } t > 0, \ z \in [a(t), a(t)+L] \tag{2.1}$$

$$(HV)_t + \left(HV^2 + \frac{1}{2}gH^2\right)_z = \mu(HV_z)_z, \text{ for } t > 0, \ z \in (a(t), a(t)+L) \tag{2.2}$$

where $H(t,z) > 0$, $V(t,z) \in \Re$ are the liquid level and the liquid velocity, respectively, at time $t \geq 0$ and position $z \in [a(t), a(t)+L]$, while $g, \mu > 0$ (constants) are the acceleration of gravity and the kinematic viscosity of the liquid, respectively.

The liquid velocities at the walls of the tank coincide with the tank velocity, i.e., we have:

$$V(t, a(t)) = V(t, a(t)+L) = w(t), \text{ for } t \geq 0 \tag{2.3}$$

where $w(t) = \dot{a}(t)$ is the velocity of the tank at time $t \geq 0$. By applying Newton's second law for the tank we get

$$\ddot{a}(t) = -f(t), \text{ for } t > 0 \tag{2.4}$$

where $-f(t)$, the control input to the problem, is the acceleration of the tank at time $t \geq 0$. The conditions for avoiding the liquid spilling out of the tank are:

$$\begin{aligned} H(t, a(t)) &< H_{max} \\ H(t, a(t)+L) &< H_{max} \end{aligned} \tag{2.5}$$

where $H_{max} > 0$ is the height of the tank walls. Applying the transformation

$$\begin{aligned} v(t,x) &= V(t, a(t)+x) - w(t) \\ h(t,x) &= H(t, a(t)+x) \\ \xi(t) &= a(t) - a^* \end{aligned} \tag{2.6}$$

where $a^* \in \Re$ is the specified position (a constant) in which we want the left wall of the tank to be placed, we obtain the model

$$\dot{\xi} = w \ , \ \dot{w} = -f, \text{ for } t \geq 0 \tag{2.7}$$

$$h_t + (hv)_x = 0, \text{ for } t > 0, \ x \in [0, L] \tag{2.8}$$

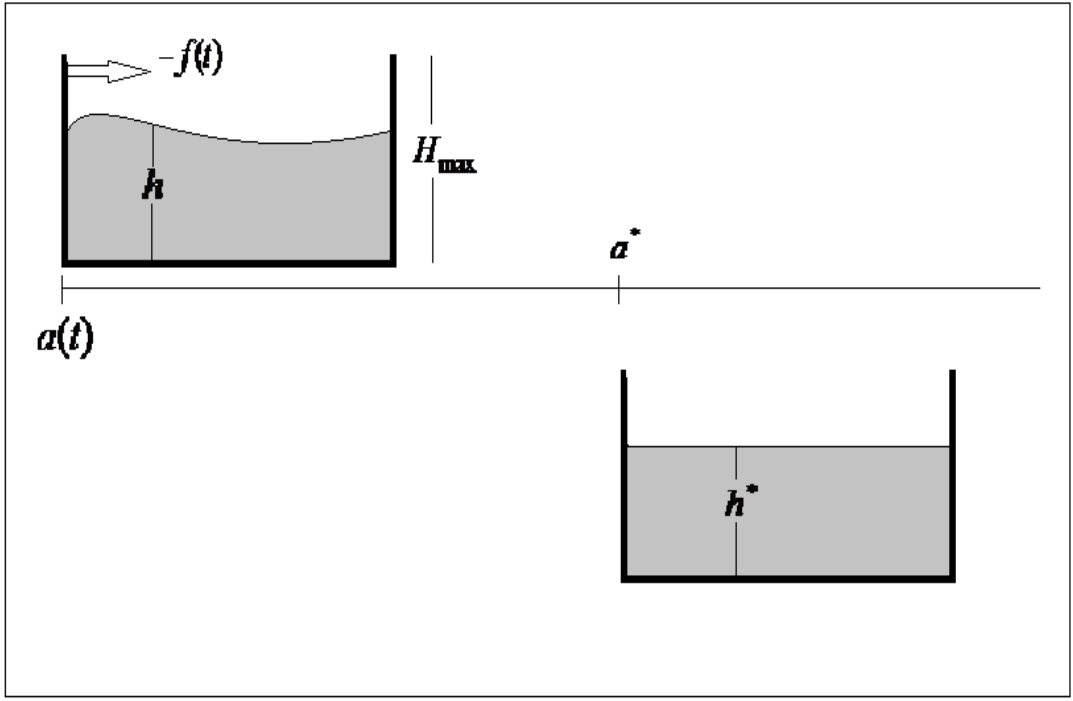

**Figure 1:** The control problem.

$$(hv)_t + \left(hv^2 + \frac{1}{2}gh^2\right)_x = \mu(hv_x)_x + hf, \text{ for } t>0,\ x\in(0,L) \quad (2.9)$$

$$v(t,0) = v(t,L) = 0, \text{ for } t \geq 0 \quad (2.10)$$

where the control input $f$ appears additively in the second equation of (2.7) and multiplicatively in (2.9). Moreover, the conditions (2.5) for avoiding the liquid spilling out of the tank become:

$$\max(h(t,0), h(t,L)) < H_{\max}, \text{ for } t \geq 0 \text{ --- Condition for no spilling out} \quad (2.11)$$

We consider classical solutions for the PDE-ODE system (2.7)-(2.10), i.e., we consider functions $\xi \in C^1(\Re_+) \cap C^2((0,+\infty))$, $w \in C^0(\Re_+) \cap C^1((0,+\infty))$, $v \in C^0([0,+\infty)\times[0,L]) \cap C^1((0,+\infty)\times[0,L])$, $h \in C^1([0,+\infty)\times[0,L];(0,+\infty)) \cap C^2((0,+\infty)\times(0,L))$, with $v[t] \in C^2((0,L))$ for each $t>0$ that satisfy equations (2.7)-(2.10) for a given input $f \in C^0(\Re_+)$.

Using (2.8) and (2.10), we can prove that for every solution of (2.7)-(2.10) it holds that $\frac{d}{dt}\left(\int_0^L h(t,x)dx\right) = 0$ for all $t>0$. Hence, the total mass of the liquid $m>0$ is constant. Therefore, without loss of generality, we assume that every solution of (2.7)-(2.10) satisfies the equation

$$\int_0^L h(t,x)dx \equiv m \quad (2.12)$$

The open-loop system (2.7)-(2.10), (2.12), i.e., system (2.7)-(2.10), (2.12) with $f(t) \equiv 0$, allows a continuum of equilibrium points, namely the points

$$h(x) \equiv h^*,\ v(x) \equiv 0, \text{ for } x\in[0,L] \quad (2.13)$$

$$\xi \in \Re,\ w = 0 \quad (2.14)$$

where $h^* = m/L$. We assume that the equilibrium points satisfy the condition for no spilling out (2.11), i.e., $h^* < H_{max}$. Due to the existence of a continuum of equilibrium points for the open-loop system given by (2.13), (2.14), the desired equilibrium point is not asymptotically stable for the open-loop system. In addition to that, experience shows that it is possible to find smooth initial conditions $(\xi(0), w(0), h[0], v[0]) \in \Re \times \Re \times C^2([0,L];(0,+\infty)) \times C^1([0,L])$ with $(v[0])(0) = (v[0])(L) = 0$ for which we cannot avoid liquid spilling out of the tank-no matter what the applied input $f$ is. Thus, the described control problem is far from trivial.

In contrast with [21], here we assume that we do not measure the whole state vector $(\xi(t), w(t), h[t], v[t])$ but we measure only the tank position $\xi(t)$ and the liquid level at the tank walls $h(t,0)$ and $h(t,L)$. In other words, the measured output is given by the equation:

$$y(t) = \begin{bmatrix} \xi(t) \\ h(t,0) \\ h(t,L) \end{bmatrix} \in \Re^3 \qquad (2.15)$$

Our objective is to design a finite-dimensional, dynamic, output feedback law of the form

$$\begin{aligned} \dot{s}(t) &= F_1(s(t), y(t)), \ s(t) \in \Re^k \\ f(t) &= F_2(s(t), y(t)) \end{aligned}, \text{ for } t > 0, \qquad (2.16)$$

that leads exponentially the solution of the closed-loop system (2.7)-(2.10) with (2.16) to the equilibrium point with $\xi = 0$. Moreover, we additionally require that the "spill-free condition" (2.11) holds for every $t \geq 0$.

## 3. Construction of the Dynamic Output Feedback Law

### *3.1. A Nonlinear State Feedback Law*

The static nonlinear state feedback law

$$f(t) = -\sigma\left(2\int_0^L h(t,x)v(t,x)dx + \mu(h(t,L) - h(t,0)) - q(w(t) + k\xi(t))\right) \qquad (3.1)$$

was proposed in [21] for the stabilization of the equilibrium point with $\xi = 0$. Here $\sigma, k, q > 0$ are constants that must satisfy the conditions

$$k < q\theta \frac{G^{-1}(-cr)}{b + G^{-1}(-cr)} \qquad (3.2)$$

where

$$c := \frac{1}{\mu\sqrt{g}}, \ \theta := \frac{\sigma g}{g + \mu\sigma L}, \ b := \frac{4mL^2 H_{max}}{\mu\pi^2}\theta, \qquad (3.3)$$

and $G^{-1}: \left(-\frac{4}{3}h^*\sqrt{h^*}, +\infty\right) \to (0, +\infty)$ is the inverse function of the increasing $C^1$ function $G: (0, +\infty) \to \left(-\frac{4}{3}h^*\sqrt{h^*}, +\infty\right)$ defined by means of the formula

$$G(h) := \frac{2}{3}\operatorname{sgn}(h-h^*)\left(h\sqrt{h} - 3h^*\sqrt{h} + 2h^*\sqrt{h^*}\right) \qquad (3.4)$$

and $r$ is a constant that satisfies

$$0 \le r < R := \frac{2\mu\sqrt{g}}{3}\left(2h^*\sqrt{h^*} + \sqrt{H_{\max}}\min\left(H_{\max} - 3h^*, 0\right)\right) \qquad (3.5)$$

The design of the feedback law (3.1) was based on the Control Lyapunov Functional (CLF) $V : S \to \Re_+$ defined by

$$V(\xi, w, h, v) := W(h,v) + E(h,v) + \frac{qk^2}{2}\xi^2 + \frac{q}{2}(w + k\xi)^2 \qquad (3.6)$$

where

$$E(h,v) := \frac{1}{2}\int_0^L h(x)v^2(x)\,dx + \frac{1}{2}g\int_0^L \left(h(x) - h^*\right)^2 dx \qquad (3.7)$$

$$W(h,v) := \frac{1}{2}\int_0^L h^{-1}(x)\left(h(x)v(x) + \mu h_x(x)\right)^2 dx + \frac{1}{2}g\int_0^L \left(h(x) - h^*\right)^2 dx \qquad (3.8)$$

and $S \subset \Re^2 \times \left(C^0([0,L])\right)^2$ is the set defined by

$$(\xi, w, h, v) \in S \Leftrightarrow \begin{cases} h \in C^0\left([0,L];(0,+\infty)\right) \cap H^1(0,L) \\ v \in C^0([0,L]) \\ \int_0^L h(x)dx = m \\ (\xi, w) \in \Re^2, v(0) = v(L) = 0 \end{cases} \qquad (3.9)$$

The feedback law (3.1) guarantees (Theorem 1 in [21]) that for each $r \in [0, R)$ there exist constants $M, \varphi > 0$ with the following property:

**(Q)** *Every classical solution of the PDE-ODE system (2.7)-(2.10), (2.12) and (3.1) with $V(\xi(0), w(0), h[0], v[0]) \le r$, satisfies $(\xi(t), w(t), h[t], v[t]) \in X$ and the following estimate for all $t \ge 0$:*

$$\left\|(\xi(t), w(t), h[t] - h^*\chi_{[0,L]}, v[t])\right\|_X \le M\exp(-\varphi t)\left\|(\xi(0), w(0), h[0] - h^*\chi_{[0,L]}, v[0])\right\|_X \qquad (3.10)$$

Here $X \subset \Re^2 \times H^1(0,L) \times L^2(0,L)$ is the metric space

$$X := \{ (\xi, w, h, v) \in S : V(\xi, w, h, v) < R \} \qquad (3.11)$$

with metric induced by the norm of the underlying normed linear space $\Re^2 \times H^1(0,L) \times L^2(0,L)$, i.e., we have for all $(\xi, w, h, v) \in X$

$$\left\|(\xi, w, h, v)\right\|_X = \left(\xi^2 + w^2 + \|h\|_2^2 + \|h'\|_2^2 + \|v\|_2^2\right)^{1/2} \qquad (3.12)$$

Finally, the "spill-free condition" (2.11) holds for every classical solution of the PDE-ODE system (2.7)-(2.10), (2.12) with (3.1) and $V(\xi(0), w(0), h[0], v[0]) \le r$, since Lemma 1 in [21] implies that for all $(\xi, w, h, v) \in X$ it holds that

$$0 < G^{-1}\left(-cV(\xi, w, h, v)\right) \leq h(x) \leq G^{-1}\left(cV(\xi, w, h, v)\right) < H_{\max}, \text{ for all } x \in [0, L] \quad (3.13)$$

The above inequality in conjunction with definitions (3.11), (3.3), (3.5) guarantees the fact that the "spill-free condition" (2.11) holds for every classical solution of the PDE-ODE system (2.7)-(2.10), (2.12) with (3.1) and $V(\xi(0), w(0), h[0], v[0]) \leq r$.

### 3.2. Ideas for the Construction of Dynamic Output Feedback Law

The main idea for the construction of a dynamic feedback law is the replacement of the quantities $\int_0^L h(t,x)v(t,x)dx$ (the liquid momentum) and $w(t)$ (the tank velocity) -which are not measured- in the feedback law (3.1) by appropriate estimates. For the tank velocity $w(t)$ an estimate $\hat{w}(t)$ can be obtained either by means of the finite-dimensional full-order observer

$$\frac{d\hat{\xi}}{dt}(t) = \hat{w}(t) - 2\gamma\left(\hat{\xi}(t) - \xi(t)\right), \quad \frac{d\hat{w}}{dt}(t) = -f(t) - \left(\frac{1}{2} + \gamma^2\right)\left(\hat{\xi}(t) - \xi(t)\right), \text{ for } t \geq 0 \quad (3.14)$$

where $\gamma > 0$ is a constant or by means of the reduced-order observer

$$\dot{\zeta}(t) = -\gamma \zeta(t) - \gamma^2 \xi(t) - f(t), \text{ for } t \geq 0 \quad (3.15)$$

where $\gamma > 0$ is a constant and the velocity estimate is given by

$$\hat{w}(t) = \zeta(t) + \gamma \xi(t) = \frac{\gamma s}{s + \gamma}[\xi] - \frac{1}{s + \gamma}[f]$$

For the liquid momentum we need to follow a different approach. If we define

$$\bar{z}(t) := \int_0^L h(t,x)v(t,x)dx + \mu\left(h(t,L) - h(t,0)\right) \quad (3.16)$$

then it follows from (2.8), (2.9), (2.10), (2.12) and definition (3.16) that the following differential equation holds for every classical solution of (2.7)-(2.10), (2.12) and all $t > 0$:

$$\frac{d\bar{z}}{dt}(t) = mf(t) - \frac{1}{2}g\left(h^2(t,L) - h^2(t,0)\right) \quad (3.17)$$

Consequently, the use of the filter

$$\dot{z}(t) = mf(t) - \frac{g}{2}\left(h^2(t,L) - h^2(t,0)\right) - \beta\left(z(t) - \mu\left(h(t,L) - h(t,0)\right)\right), \text{ for } t > 0 \quad (3.18)$$

where $\beta > 0$ is a constant, along with (3.17), guarantees the following differential equation for every classical solution of (2.7)-(2.10), (2.12) and all $t > 0$:

$$\frac{d}{dt}\left(z(t) - \bar{z}(t)\right) = -\beta\left(z(t) - \bar{z}(t)\right) - \beta\int_0^L h(t,x)v(t,x)dx \quad (3.19)$$

We can now estimate the liquid momentum $p = \int_0^L h(t,x)v(t,x)dx$ by:

- either the quantity

$$\hat{p} = z(t) - \mu(h(t,L) - h(t,0)) = \frac{1}{s+\beta}\left[mf - \frac{g}{2}(h^2(L) - h^2(0))\right] - \frac{s}{s+\beta}\left[\mu(h(L) - h(0))\right],$$

  which is going to give an accurate estimate when $|z(t) - \bar{z}(t)|$ is small,

- or by zero, which is going to give an accurate estimate when $\left|\int_0^L h(t,x)v(t,x)dx\right|$ is small.

Taking the convex combination of both estimates of the liquid momentum we get the estimate $\lambda z(t) - \lambda \mu(h(t,L) - h(t,0))$, where $\lambda \in [0,1]$ is a constant.

The reader should notice that when $\lambda = 0$ then the filter (3.18) is unutilized. Thus we may obtain four different dynamic output feedback laws with:

(a) three ODEs when the full-order observer (3.14) is used and $\lambda \in (0,1]$,

(b) two ODEs when the reduced-order observer (3.15) is used and $\lambda \in (0,1]$,

(c) two ODEs when the full-order observer (3.14) is used and $\lambda = 0$, and

(d) one ODE when the reduced-order observer (3.15) is used and $\lambda = 0$.

## 4. Main Results

Now we are in a position to state the main results of the present work.

**Theorem 1 (Output Feedback with Full-Order Observer):** *For every $r \in [0, R)$, where $R > 0$ is the constant defined in (3.5), and for every $\sigma, q, k, \gamma, \beta > 0$, $\lambda \in [0,1]$ with*

$$\frac{2g}{\mu L} > \sigma > \frac{8k}{q}, \quad \frac{2\pi^2 \mu G^{-1}(-cr)}{mL^2 G^{-1}(cr)} > \beta + 4\sigma > 8\sigma \qquad (4.1)$$

*where $c > 0$ is defined by (3.3), there exist constants $M, \varphi > 0$ with the following property:*

**(P)** *Every classical solution of the PDE-ODE system (2.7)-(2.10), (2.12), (3.14), (3.18) and*

$$f = -\sigma\left(2\lambda z + (1-2\lambda)\mu(h(L) - h(0)) - q(\hat{w} + k\xi) + q\gamma(\hat{\xi} - \xi)\right), \text{ for } t > 0 \qquad (4.2)$$

*with $\Phi(\xi(0), w(0), h[0], v[0], \hat{\xi}(0), \hat{w}(0), z(0)) \le r$, where $\Phi: X \times \mathbb{R}^3 \to \mathbb{R}_+$ is the functional*

$$\Phi(\xi, w, h, v, \hat{\xi}, \hat{w}, z) := V(\xi, w, h, v) + \frac{\sigma q^2}{2\gamma}(\hat{\xi} - \xi)^2$$
$$+ \frac{\sigma q^2}{\gamma}\left(\hat{w} - w - \gamma(\hat{\xi} - \xi)\right)^2 + \frac{\lambda}{2}\left(z - \int_0^L h(x)v(x)dx - \mu(h(L) - h(0))\right)^2 \qquad (4.3)$$

*satisfies $(\xi(t), w(t), h[t], v[t]) \in X$ and the following estimate for all $t \ge 0$:*

$$\left\|(\xi(t),w(t),h[t]-h^*\chi_{[0,L]},v[t])\right\|_X + \left|\hat{\xi}(t)-\xi(t)\right|$$

$$+\left|\hat{w}(t)-w(t)\right|+\sqrt{\lambda}\left|z(t)-\int_0^L h(t,x)v(t,x)dx-\mu\big(h(t,L)-h(t,0)\big)\right|$$

$$\leq M\exp(-\varphi t)\left(\left\|(\xi(0),w(0),h[0]-h^*\chi_{[0,L]},v[0])\right\|_X+\left|\hat{\xi}(0)-\xi(0)\right|\right) \quad (4.4)$$

$$+M\exp(-\varphi t)\left(\left|\hat{w}(0)-w(0)\right|+\sqrt{\lambda}\left|z(0)-\int_0^L h(0,x)v(0,x)dx-\mu\big(h(0,L)-h(0,0)\big)\right|\right)$$

**Theorem 2 (Output Feedback with Reduced-Order Observer):** *For every $r \in [0,R)$, where $R > 0$ is the constant defined in (3.5), and for every $\sigma, q, k, \gamma, \beta > 0$, $\lambda \in [0,1]$ satisfying (4.1) there exist constants $M, \varphi > 0$ with the following property:*

**(P')** *Every classical solution of the PDE-ODE system (2.7)-(2.10), (2.12), (3.15), (3.18) and*

$$f = -\sigma\big(2\lambda z + (1-2\lambda)\mu(h(L)-h(0)) - q(\zeta + (\gamma+k)\xi)\big), \text{ for } t > 0 \quad (4.5)$$

*with $\Psi(\xi(0), w(0), h[0], v[0], \zeta(0), z(0)) \leq r$, where $\Psi: X \times \mathfrak{R}^2 \to \mathfrak{R}_+$ is the functional*

$$\Psi(\xi,w,h,v,\zeta,z) := V(\xi,w,h,v) + \frac{\sigma q^2}{\gamma}(\zeta - w + \gamma\xi)^2$$

$$+\frac{\lambda}{2}\left(z-\int_0^L h(x)v(x)dx-\mu\big(h(L)-h(0)\big)\right)^2 \quad (4.6)$$

*satisfies $(\xi(t), w(t), h[t], v[t]) \in X$ and the following estimate for all $t \geq 0$:*

$$\left\|(\xi(t),w(t),h[t]-h^*\chi_{[0,L]},v[t])\right\|_X + \left|\zeta(t)-w(t)+\gamma\xi(t)\right|$$

$$+\sqrt{\lambda}\left|z(t)-\int_0^L h(t,x)v(t,x)dx-\mu\big(h(t,L)-h(t,0)\big)\right|$$

$$\leq M\exp(-\varphi t)\left(\left\|(\xi(0),w(0),h[0]-h^*\chi_{[0,L]},v[0])\right\|_X\right. \quad (4.7)$$

$$+\left|\zeta(0)-w(0)+\gamma\xi(0)\right|$$

$$\left.+\sqrt{\lambda}\left|z(0)-\int_0^L h(0,x)v(0,x)dx-\mu\big(h(0,L)-h(0,0)\big)\right|\right)$$

**Remarks: (a)** Both Theorem 1 and Theorem 2 guarantee exponential convergence to the desired equilibrium point. Moreover, both estimates (4.4) and (4.7) are independent of $z$ when $\lambda = 0$. This is expected since the corresponding closed-loop systems are independent of $z$ when $\lambda = 0$.
**(b)** It should be noticed that Theorem 1 in [21] had no condition whatsoever for the controller gain $\sigma > 0$ in the feedback laws (4.2) and (4.5). However, Theorem 1 and Theorem 2 in the present work require the validity of inequality (4.1) which shows that $\sigma > 0$ should be sufficiently large but cannot be arbitrarily large. The reason for the existence of an upper bound for $\sigma > 0$ is the existence of estimation errors in the present work. Since we cannot implement the state feedback law (3.1) but instead we try to approximate the feedback law (3.1) by using appropriate estimates, the use of a large controller gain $\sigma > 0$ will also magnify the effect of the estimation errors to the closed-loop systems.

**(c)** On the other hand, neither Theorem 1 nor Theorem 2 require any upper bound for the parameter $\gamma > 0$ that determines the convergence rate to zero of the observer error for the full observer (3.14) and the reduced observer (3.15). However, a very large value for the observer gain $\gamma > 0$ should not be used in practice due to the existence of possible measurement noise.

**(d)** The gain $\beta > 0$ of the filter (3.18) has to be greater than $4\sigma$ and less than $\dfrac{2\pi^2 \mu G^{-1}(-cr)}{mL^2 G^{-1}(cr)} - 4\sigma$ (recall (4.1)). Indeed, the use of a large filter gain $\beta > 0$ will magnify the effect of the estimation errors to both closed-loop systems and this is the reason that the gain $\beta$ has to be less than $\dfrac{2\pi^2 \mu G^{-1}(-cr)}{mL^2 G^{-1}(cr)} - 4\sigma$.

## 5. Can a Robot Move a Glass of Water Without Spilling Out Water and Without Measuring the Water Momentum?

We next consider the problem of the movement of a glass of water by means of a robotic arm, from near-rest to near-rest in finite (but not prespecified) time. The glass of water starting from an almost at rest state (both the glass and the water in the glass), is formulated as

$$\left\|(0, w(0), h[0] - h^* \chi_{[0,L]}, v[0])\right\|_X \le \varepsilon \tag{5.1}$$

where $\varepsilon > 0$ is a small number (tolerance) and $\|\cdot\|_X$ is the norm defined by (3.12). The robotic arm should move the glass to a pre-specified position without spilling out water of the glass and having the water in the glass almost still at the final time, i.e., at the final time $T > 0$ we must have

$$\left\|(\xi(T), w(T), h[T] - h^* \chi_{[0,L]}, v[T])\right\|_X \le \varepsilon \tag{5.2}$$

In other words, we require the spill-free and slosh-free motion of the glass. The initial condition is $(\xi(0), w(0), h[0], v[0]) \in S$, where $\xi(0) \ne 0$ (recall the definition of $S$ (3.9)).

The problem was solved in [21] for a sufficiently small tolerance $\varepsilon > 0$ by using the feedback law (3.1). However, the feedback law (3.1) requires the measurement both of the momentum of the water in the glass and the glass velocity. Both these two measurements are not easily obtained in practice and thus the problem of performing of spill-free and slosh-free transfer of the glass of water without measuring the momentum of the water in the glass and without measuring the glass velocity is important.

Here, we solve the problem by applying Theorem 2 with $\lambda = 0$ when the tolerance $\varepsilon > 0$ is sufficiently small (but also Theorem 1 or different values of $\lambda \in [0,1]$ can be applied). More specifically, we require that $\varepsilon > 0$ is small so that

$$\varepsilon^2 \max\left(\mu^2 \left(h^* - \varepsilon\sqrt{L}\right)^{-1}, g, \frac{3H_{\max}}{2}\right) < R \text{ and } \varepsilon < \frac{\min\left(h^*, H_{\max} - h^*\right)}{\sqrt{L}} \tag{5.3}$$

where $R := \dfrac{2\mu\sqrt{g}}{3}\left(2h^*\sqrt{h^*} + \sqrt{H_{\max}}\min\left(H_{\max} - 3h^*, 0\right)\right)$.

Inequalities (5.3) are exactly the same conditions for tolerance as in the case of the state feedback law (3.1) that was studied in [21]. If inequalities (5.3) hold then we can follow the following algorithm:

**Step 1:** Pick numbers $r \in (0, R)$, $q > 0$ with $q \leq \max\left(\mu^2\left(h^* - \varepsilon\sqrt{L}\right)^{-1}, g, \frac{3H_{max}}{2}\right)$ and $\varepsilon^2 \max\left(\mu^2\left(h^* - \varepsilon\sqrt{L}\right)^{-1}, g, \frac{3H_{max}}{2}\right) < r$ (this is possible due to (5.3)).

**Step 2:** Select $\sigma, k, \gamma > 0$ so that the inequalities $\frac{2g}{\mu L} > \sigma > \frac{8k}{q}$, $\frac{2\pi^2 \mu G^{-1}(-cr)}{mL^2 G^{-1}(cr)} > 8\sigma$ hold and so that

$$\gamma \geq \frac{2\sigma q^2 \varepsilon^2}{r - \max\left(\mu^2\left(h^* - \varepsilon\sqrt{L}\right)^{-1}, g, \frac{3H_{max}}{2}\right)\varepsilon^2} \quad \text{and} \quad k \leq \sqrt{\frac{1}{3q}} \frac{\sqrt{r - \varepsilon^2 \max\left(\mu^2\left(h^* - \varepsilon\sqrt{L}\right)^{-1}, g, \frac{3H_{max}}{2}\right)}}{|\xi(0)|}$$

(5.4)

Inequalities (5.4) guarantee the inequality $\Psi(\xi(0), w(0), h[0], v[0], \zeta(0), z(0)) \leq r$ for $\zeta(0) = -\gamma\xi(0)$ and for arbitrary $z(0) \in \Re$ (since $\lambda = 0$ the filter (3.18) is irrelevant). Indeed, this fact is a direct consequence of (3.12), (4.6), (5.1), (5.3), (5.4) and Proposition 1 in [21]. More specifically, since $\varepsilon < \frac{\min(h^*, H_{max} - h^*)}{\sqrt{L}}$ and $(\xi(0), w(0), h[0], v[0]) \in S$ satisfies (5.1), Proposition 1 in [21] implies the inequality

$$V(\xi(0), w(0), h[0], v[0])$$
$$\leq \max\left(\mu^2\left(h^* - \varepsilon\sqrt{L}\right)^{-1}, g, \frac{3H_{max}}{2}, q\right) \left\|(0, w(0), h[0] - h^*\chi_{[0,L]}, v[0])\right\|_X^2 + \frac{3qk^2}{2}\xi^2(0)$$

Since $q \leq \max\left(\mu^2\left(h^* - \varepsilon\sqrt{L}\right)^{-1}, g, \frac{3H_{max}}{2}\right)$, the above inequality (5.1) and definition (4.6) gives the estimate (recall that $\lambda = 0$):

$$\Psi(\xi(0), w(0), h[0], v[0], \zeta(0), z(0))$$
$$= V(\xi(0), w(0), h[0], v[0]) + \frac{\sigma q^2}{\gamma}(\zeta(0) - w(0) + \gamma\xi(0))^2$$
$$\leq \max\left(\mu^2\left(h^* - \varepsilon\sqrt{L}\right)^{-1}, g, \frac{3H_{max}}{2}\right)\varepsilon^2 + \frac{3qk^2}{2}\xi^2(0) + \frac{\sigma q^2}{\gamma}(\zeta(0) - w(0) + \gamma\xi(0))^2$$

The above inequality in conjunction with (5.4) and the fact that $|w(0)| \leq \left\|(0, w(0), h[0] - h^*\chi_{[0,L]}, v[0])\right\|_X \leq \varepsilon$ (a consequence of (3.12) and (5.1)) guarantees that $\Psi(\xi(0), w(0), h[0], v[0], \zeta(0), z(0)) \leq r$ for $\zeta(0) = -\gamma\xi(0)$ and for arbitrary $z(0) \in \Re$.

Notice that since $\frac{2g}{\mu L} > \sigma > \frac{8k}{q}$, $\frac{2\pi^2 \mu G^{-1}(-cr)}{mL^2 G^{-1}(cr)} > 8\sigma$ inequalities (4.1) hold for appropriate $\beta > 0$ but notice that $\beta > 0$ in this case is irrelevant due to the fact that $\lambda = 0$ and the filter (3.18) is not used. Consequently, we are in a position to apply Theorem 2. Let $M, \varphi > 0$ be the constants involved in (4.7) and correspond to the selected parameters $r \in (0, R)$, $\lambda = 0$ and $\sigma, k, \gamma, \beta > 0$.

**Step 3:** Set $T = \frac{1}{\varphi} \ln\left( \frac{M|\xi(0)| + 2M\varepsilon}{\varepsilon} \right)$ and apply the feedback law (3.15), (4.5) with $\zeta(0) = -\gamma\xi(0)$, $\lambda = 0$ for $t \in (0,T]$. Inequalities (3.12), (4.7), (5.1) imply that estimate (5.2) holds. Indeed, we obtain from (4.7) for $\lambda = 0$:

$$\|(\xi(T), w(T), h[T] - h^*\chi_{[0,L]}, v[T])\|_X$$
$$\leq M \exp(-\varphi T)\left( \|(\xi(0), w(0), h[0] - h^*\chi_{[0,L]}, v[0])\|_X + |\zeta(0) - w(0) + \gamma\xi(0)| \right)$$

The above inequality and the fact that $\zeta(0) = -\gamma\xi(0)$ in conjunction with definition (3.12) (which directly implies that $\|(\xi(0), w(0), h[0] - h^*\chi_{[0,L]}, v[0])\|_X \leq |\xi(0)| + \|(0, w(0), h[0] - h^*\chi_{[0,L]}, v[0])\|_X$) gives the estimate

$$\|(\xi(T), w(T), h[T] - h^*\chi_{[0,L]}, v[T])\|_X$$
$$\leq M \exp(-\varphi T)\left( \|(0, w(0), h[0] - h^*\chi_{[0,L]}, v[0])\|_X + |\xi(0)| + |w(0)| \right)$$

Finally, inequality (5.1) and the fact that $|w(0)| \leq \|(0, w(0), h[0] - h^*\chi_{[0,L]}, v[0])\|_X \leq \varepsilon$ (a consequence of (3.12) and (5.1)) imply the estimate

$$\|(\xi(T), w(T), h[T] - h^*\chi_{[0,L]}, v[T])\|_X \leq M \exp(-\varphi T)\left( |\xi(0)| + 2\varepsilon \right)$$

Estimate (5.2) is a direct consequence of the above estimate and the fact that $T = \frac{1}{\varphi} \ln\left( \frac{M|\xi(0)| + 2M\varepsilon}{\varepsilon} \right)$.

The application of the above algorithm can guarantee that the robotic arm will move the glass to the pre-specified position *without spilling out water of the glass, without residual end point sloshing, without measuring the glass velocity and without measuring the momentum of water in the glass - no matter how small the difference* $H_{\max} - h^* > 0$ *is*. Of course, the (minimization of the) final time $T > 0$ depends on the (appropriate) selection of the parameters $r \in (0, R)$ and $\sigma, q, k, \gamma > 0$.

## 6. Numerical Examples

Applying appropriate scaling and using dimensionless variables in the viscous Saint-Venant model (2.7)-(2.10), (2.12), we may assume (without any loss of generality) next that $L = g = h^* = 1$.

Using finite differences for (2.7), (2.8), (2.9) and (2.10) with time step $\delta t > 0$, spatial discretization step $\delta x = 1/n$, where $n \geq 2$ is an integer, and using the notation $h_i, v_i$ for the values of the numerical approximations of the liquid level and velocity, respectively, at position $x = i\delta x$ with $i = 0, \ldots, n$, we obtain the following discrete-time system (numerical scheme):

$$h_0^+ = h_0 \exp\left(\frac{\delta t}{2\delta x}(v_2 - 4v_1)\right),$$

$$h_i^+ = h_i \exp\left(\frac{\delta t}{2\delta x}(v_{i-1} - v_{i+1}) + \frac{\delta t}{2\delta x} v_i \ln\left(\frac{h_{i-1}}{h_{i+1}}\right)\right) \text{ for } i = 1,\ldots,n-1, \qquad (6.1)$$

$$h_n^+ = h_n \exp\left(\frac{\delta t}{2\delta x}(4v_{n-1} - v_{n-2})\right)$$

and

$$v_i^+ = \left(1 - \frac{2\delta t \mu}{(\delta x)^2}\right) v_i + \frac{\delta t \mu}{(\delta x)^2}(v_{i+1} + v_{i-1}) + \frac{\delta t \mu}{4(\delta x)^2}(v_{i+1} - v_{i-1}) \ln\left(\frac{h_{i+1}}{h_{i-1}}\right)$$

$$- \frac{\delta t}{2\delta x}\left(v_i(v_{i+1} - v_{i-1}) + h_{i+1} - h_{i-1}\right) + f \delta t \qquad (6.2)$$

for $i = 1,\ldots,n-1$ with

$$v_0^+ = v_n^+ = 0 \qquad (6.3)$$

In the case of the output feedback with full-order observer the above model is combined with the control law

$$f = -\sigma\left(2\lambda z + (1 - 2\lambda)\mu(h_n - h_0) - q(\hat{w} + k\hat{\xi}) + q\gamma(\hat{\xi} - \xi)\right) \qquad (6.4)$$

where the velocity of the tank $w$, the position of the tank $\xi$, the velocity estimate $\hat{w}$ and the position estimate $\hat{\xi}$ at the next time instant are given by the following equations

$$\xi^+ = \xi + \delta t\, w - \frac{(\delta t)^2}{2} f$$

$$w^+ = w - \delta t\, f$$

$$\hat{\xi}^+ = \xi^+ + \exp(-\gamma\, \delta t)\left(\cos\left(\frac{\delta t}{\sqrt{2}}\right)(\hat{\xi} - \xi) + \sqrt{2}\sin\left(\frac{\delta t}{\sqrt{2}}\right)(\hat{w} - w - \gamma(\hat{\xi} - \xi))\right), \qquad (6.5)$$

$$\hat{w}^+ = w^+ + \exp(-\gamma\, \delta t)\left(\cos\left(\frac{\delta t}{\sqrt{2}}\right)(\hat{w} - w) + \gamma\sqrt{2}\sin\left(\frac{\delta t}{\sqrt{2}}\right)\left(\hat{w} - w - \left(\gamma + \frac{1}{2\gamma}\right)(\hat{\xi} - \xi)\right)\right)$$

In the case of the reduced-order observer the model is combined with the control law

$$f = -\sigma\left(2\lambda z + (1 - 2\lambda)\mu(h_n - h_0) - q(\zeta + (\gamma + k)\xi)\right) \qquad (6.6)$$

where the velocity of the tank $w$, the position of the tank $\xi$ and the observer state $\zeta$ at the next time instant are given by the following scheme

$$\xi^+ = \xi + \delta t\, w - \frac{(\delta t)^2}{2} f$$

$$w^+ = w - \delta t\, f \qquad (6.7)$$

$$\zeta^+ = w^+ - \gamma\xi^+ + \exp(-\gamma\, \delta t)(\zeta - w + \gamma\xi)$$

In both cases the filter state $z$ at the next time instant is given by the following scheme

$$z^+ = z + \delta t \left( f - \frac{1}{2}\left(h_n^2 - h_0^2\right) - \beta\left(z - \mu\left(h_n - h_0\right)\right)\right) \quad (6.8)$$

The discretization schemes (6.1)-(6.5), (6.8) for the viscous Saint-Venant system under the output feedback controller with full-order observer and (6.1)-(6.3), (6.6)-(6.8) for the viscous Saint-Venant system under the output feedback controller with reduced-order observer are non-standard, explicit, finite-difference schemes and they guarantee positivity of the liquid level (i.e., if $h_i > 0$ for $i = 0,\ldots,n$ then $h_i^+ > 0$ for $i = 0,\ldots,n$). Moreover, the discretization schemes (6.1)-(6.5), (6.8) and (6.1)-(6.3), (6.6)-(6.8) are accurate of order (2,1) (see Definition 2.3.3 on page 57 in the book [43]). We have not been able to find similar discretization schemes in the literature and therefore, we have included the following results for completeness (their proofs are given in the following section). By the term "smooth solution" we mean a solution with derivatives of sufficiently large order.

**Theorem 3 (Numerical Discretization of the Viscous Saint-Venant System Under the Output Feedback with Full-Order Observer):** *Let constants $\sigma, q, k, \gamma, \beta > 0$, $\lambda \in [0,1]$ be given. Let $T > 0$ be a given constant and let $(\xi(t), w(t), h[t], v[t], \hat{\xi}(t), \hat{w}(t), z(t)) \in X \times \Re^3$ with $t \in [0,T]$ be a classical smooth solution of the PDE-ODE system (2.7)-(2.10), (2.12), (3.14), (3.18), (4.2) with $L = g = h^* = 1$. Then there exists a constant $\bar{M} > 0$ such that the following estimate holds for all sufficiently small $\delta t \in (0,T)$, $\delta x = 1/n$*

$$\left|\xi^+ - \xi(\delta t)\right| + \left|w^+ - w(\delta t)\right| + \left|z^+ - z(\delta t)\right| + \left|\hat{\xi}^+ - \hat{\xi}(\delta t)\right| + \left|\hat{w}^+ - \hat{w}(\delta t)\right|$$
$$+ \max_{i=0,\ldots,n}\left(\left|h_i^+ - h(\delta t, i\delta x)\right| + \left|v_i^+ - v(\delta t, i\delta x)\right|\right) \leq \bar{M}\,\delta t\left(\delta t + (\delta x)^2\right) \quad (6.9)$$

*where $n \geq 2$ is an integer and $h_i^+, v_i^+$ for $i = 0,\ldots,n$, $w^+, \xi^+, z^+, \hat{\xi}^+, \hat{w}^+$ are given by the difference scheme (6.1)-(6.5), (6.8) with $h_i = h(0, i\delta x)$, $v_i = v(0, i\delta x)$, for $i = 0,\ldots,n$, $w = w(0)$, $\xi = \xi(0)$, $\hat{\xi} = \hat{\xi}(0)$, $\hat{w} = \hat{w}(0)$ and $z = z(0)$.*

**Theorem 4 (Numerical Discretization of the Viscous Saint-Venant System Under the Output Feedback with Reduced-Order Observer):** *Let constants $\sigma, q, k, \gamma, \beta > 0$, $\lambda \in [0,1]$ be given. Let $T > 0$ be a given constant and let $(\xi(t), w(t), h[t], v[t], \zeta(t), z(t)) \in X \times \Re^2$ with $t \in [0,T]$ be a classical smooth solution of the PDE-ODE system (2.7)-(2.10), (2.12), (3.15), (3.18), (4.5) with $L = g = h^* = 1$. Then there exists a constant $\bar{M} > 0$ such that the following estimate holds for all sufficiently small $\delta t \in (0,T)$, $\delta x = 1/n$*

$$\left|\xi^+ - \xi(\delta t)\right| + \left|w^+ - w(\delta t)\right| + \left|z^+ - z(\delta t)\right| + \left|\zeta^+ - \zeta(\delta t)\right|$$
$$+ \max_{i=0,\ldots,n}\left(\left|h_i^+ - h(\delta t, i\delta x)\right| + \left|v_i^+ - v(\delta t, i\delta x)\right|\right) \leq \bar{M}\,\delta t\left(\delta t + (\delta x)^2\right) \quad (6.10)$$

*where $n \geq 2$ is an integer and $h_i^+, v_i^+$ for $i = 0,\ldots,n$, $w^+, \xi^+, z^+, \zeta^+$ are given by the difference scheme (6.1)-(6.3), (6.6)-(6.8) with $h_i = h(0, i\delta x)$, $v_i = v(0, i\delta x)$, for $i = 0,\ldots,n$, $w = w(0)$, $\xi = \xi(0)$, $\zeta = \zeta(0)$ and $z = z(0)$.*

We have performed many simulations with $n = 100$ for various values of the parameters $\sigma, q, k, \gamma, \beta > 0$, $\lambda \in [0,1]$. For each simulation we tested whether the inequalities

$$\bar{\Phi}^+ \leq \bar{\Phi}, \text{ for the output feedback with full observer}$$

$$\bar{\Psi}^+ \leq \bar{\Psi}, \text{ for the output feedback with reduced observer}$$

were valid for all times, where

$$\bar{\Phi} = \frac{qk^2}{2}\xi^2 + \frac{q}{2}(w+k\xi)^2 + \frac{1}{2n}\sum_{i=1}^{n-1}h_i v_i^2 + \frac{1}{n}\sum_{i=1}^{n-1}(h_i-1)^2 + \frac{1}{2n}(h_0-1)^2 + \frac{1}{2n}(h_n-1)^2$$

$$+\frac{\sigma q^2}{2\gamma}(\hat{\xi}-\xi)^2 + \frac{\sigma q^2}{\gamma}\left(\hat{w}-w-\gamma(\hat{\xi}-\xi)\right)^2 + \frac{\lambda}{2}\left(z - \frac{1}{n}\sum_{i=1}^{n-1}h_i v_i - \mu(h_n - h_0)\right)^2 \quad (6.11)$$

$$+\frac{1}{2n}\sum_{i=1}^{n-1}h_i^{-1}\left(h_i v_i + \frac{n\mu}{2}(h_{i+1}-h_{i-1})\right)^2 + \frac{n\mu^2}{16}h_0^{-1}(4h_1 - h_2 - 3h_0)^2 + \frac{n\mu^2}{16}h_n^{-1}(3h_0 + h_{n-2} - 4h_{n-1})^2$$

$$\bar{\Psi} = \frac{\sigma q^2}{\gamma}(\zeta - w + \gamma\xi)^2 + \frac{q}{2}(w+k\xi)^2 + \frac{qk^2}{2}\xi^2 + \frac{\lambda}{2}\left(z - \frac{1}{n}\sum_{i=1}^{n-1}h_i v_i - \mu(h_n-h_0)\right)^2$$

$$+\frac{1}{n}\sum_{i=1}^{n-1}(h_i-1)^2 + \frac{1}{2n}(h_0-1)^2 + \frac{1}{2n}(h_n-1)^2 + \frac{1}{2n}\sum_{i=1}^{n-1}h_i^{-1}\left(h_i v_i + \frac{n\mu}{2}(h_{i+1}-h_{i-1})\right)^2 \quad (6.12)$$

$$+\frac{n\mu^2}{16}h_0^{-1}(4h_1 - h_2 - 3h_0)^2 + \frac{n\mu^2}{16}h_n^{-1}(3h_0 + h_{n-2} - 4h_{n-1})^2 + \frac{1}{2n}\sum_{i=1}^{n-1}h_i v_i^2$$

and $\bar{\Phi}^+$, $\bar{\Psi}^+$ are the values of the above functions with $h_i, v_i$ for $i = 0,1,\ldots,n$ replaced by $h_i^+, v_i^+$, respectively, for $i = 0,1,\ldots,n$, $\xi, w, z$ replaced by $\xi^+, w^+, z^+$, $\hat{\xi}, \hat{w}$ replaced by $\hat{\xi}^+, \hat{w}^+$ for the output feedback with full observer and $\zeta$ replaced by $\zeta^+$ for the output feedback with reduced observer. Notice that $\bar{\Phi}, \bar{\Psi}$ are numerical approximations of the values of the functionals $\Phi, \Psi$ defined by (4.3), (4.6), respectively (with all integrals evaluated by means of the trapezoid rule and derivatives approximated by central differences). The proofs of Theorem 1 and Theorem 2 show that the functionals $\Phi, \Psi$ are strict Lyapunov functionals for the corresponding closed-loop systems and thus the inequality $\bar{\Phi}^+ \leq \bar{\Phi}$ or $\bar{\Psi}^+ \leq \bar{\Psi}$ must necessarily hold for each time, provided that the time step $\delta t > 0$ is appropriately selected. If there were a time for which the inequality $\bar{\Phi}^+ \leq \bar{\Phi}$ or $\bar{\Psi}^+ \leq \bar{\Psi}$ was not valid then the simulation was not accepted and we reduced the time step $\delta t > 0$.

In order to study the convergence of the state we used the following non-negative function

$$\Omega := \left(\xi^2 + w^2 + \frac{1}{n}\sum_{i=1}^{n-1}(h_i-1)^2 + \frac{1}{2n}(h_0-1)^2 + \frac{1}{2n}(h_n-1)^2 + \frac{1}{n}\sum_{i=1}^{n-1}v_i^2\right)^{1/2} \quad (6.13)$$

which is a numerical approximation (by means of the trapezoidal rule) of the functional $\left(\xi^2 + w^2 + \|h - \chi_{[0,1]}\|_2^2 + \|v\|_2^2\right)^{1/2}$. We found that by increasing the controller gain $\sigma$ and the observer gain $\gamma$ we were able to guarantee faster convergence. The controller parameters $q, k, \beta > 0$ play a less significant role. The parameter $\lambda \in [0,1]$ (the weight of the filter estimation) played a minor role and in many cases the time evolution of the state under the proposed output feedback laws with $\lambda = 0$ and $\lambda = 1$ was almost identical. Figure 2 and Figure 3 show the time evolution of the function $\Omega$ defined by (6.13) with $\delta t = 10^{-4}$, $\mu = 0.2$, $\gamma = 1$, $\sigma = 9.5$, $q = 1.2$, $k = 1.5$, $\lambda = 0$ and

$$\xi(0) = 2, \ w(0) = -2.2, \ h(0,x) = x + 0.5, \ v(0,x) = x^2 - x, \ z(0) = 0.2 \text{ and } \zeta(0) = -5$$
for the output feedback with the reduced-order observer

$$\xi(0) = 2, \ w(0) = -2.2, \ \hat{\xi}(0) = 2.3, \ \hat{w}(0) = -2.5, \ h(0,x) = x + 0.5, \ v(0,x) = x^2 - x \text{ and } z(0) = 0.2$$
for the output feedback with the full-order observer

Figure 3 shows the time evolution of the logarithms of the functions $\bar{\Phi}$ (for the case of output feedback with the full-order observer), $\bar{\Psi}$ (for the case of output feedback with the reduced-order observer) defined by (6.11), (6.12), respectively. The plots clearly indicate exponential convergence. This is expected since the proofs of Theorem 1 and Theorem 2 show exponential convergence of the functionals $\Phi, \Psi$ defined by (4.3), (4.6), respectively (recall that $\bar{\Phi}, \bar{\Psi}$ are numerical approximations of the values of the functionals $\Phi, \Psi$).

Figure 4 shows snapshots of the liquid level profile for the case of output feedback with the reduced-order observer. It is shown that the liquid level approaches quickly the equilibrium although sloshing during an initial transient is apparent.

Finally, it should be noticed that simulations showed that for both output feedback schemes and for fixed controller parameters $\sigma, q, k, \gamma, \beta > 0$, $\lambda \in [0,1]$ there is an optimal value $\mu_{opt}$ of the viscosity coefficient $\mu$ for which the convergence of the function $\Omega$ defined by (6.13) becomes fastest. It is important to mention that $\mu_{opt}$ is different for each output feedback scheme and for different combinations of the controller parameters $\sigma, q, k, \gamma, \beta$ and $\lambda$.

## 7. Proofs

We first give the proof of Theorem 1.

**Proof of Theorem 1:** Let $r \in [0, R)$ be given (arbitrary) and let constants $\sigma, q, k, \gamma > 0$, $\lambda \in [0,1]$ for which (4.1) holds (but otherwise arbitrary) be given. Consider a classical solution of the PDE-ODE system (2.7)-(2.10), (2.12) with (3.14), (3.18), (4.2) that satisfies $\Phi(\xi(0), w(0), h[0], v[0], \hat{\xi}(0), \hat{w}(0), z(0)) \leq r$.

Using Lemma 2 in [21], (2.7), (3.14), (3.19) and definitions (3.16), (4.3), we obtain for all $t > 0$:

$$\begin{aligned}
\dot{\Phi} = & -\mu g \int_0^L h_x^2(x)dx - \mu \int_0^L h(x)v_x^2(x)dx - qk^3\xi^2 \\
& + qk(w+k\xi)^2 + f\left(2\bar{z} - \mu(h(L)-h(0)) - q(w+k\xi)\right) \\
& - \sigma q^2\left(\hat{\xi}-\xi\right)^2 - 2\sigma q^2\left(\hat{w}-w-\gamma\left(\hat{\xi}-\xi\right)\right)^2 \\
& - \beta\lambda\left(z-\bar{z}\right)^2 - \beta\lambda\left(z-\bar{z}\right)\int_0^L h(x)v(x)dx
\end{aligned} \qquad (7.1)$$

Equation (4.2) and definition (3.16) implies that

$$\begin{aligned}
f = & -\sigma\left(2\bar{z} - \mu(h(L)-h(0)) - q(w+k\xi)\right) - 2\sigma\lambda(z-\bar{z}) \\
& + \sigma q\left(\hat{w}-w-\gamma\left(\hat{\xi}-\xi\right)\right) + 2\sigma(1-\lambda)\int_0^L h(x)v(x)dx
\end{aligned} \qquad (7.2)$$

for all $t > 0$. Combining equations (7.1), (7.2), we obtain for all $t > 0$:

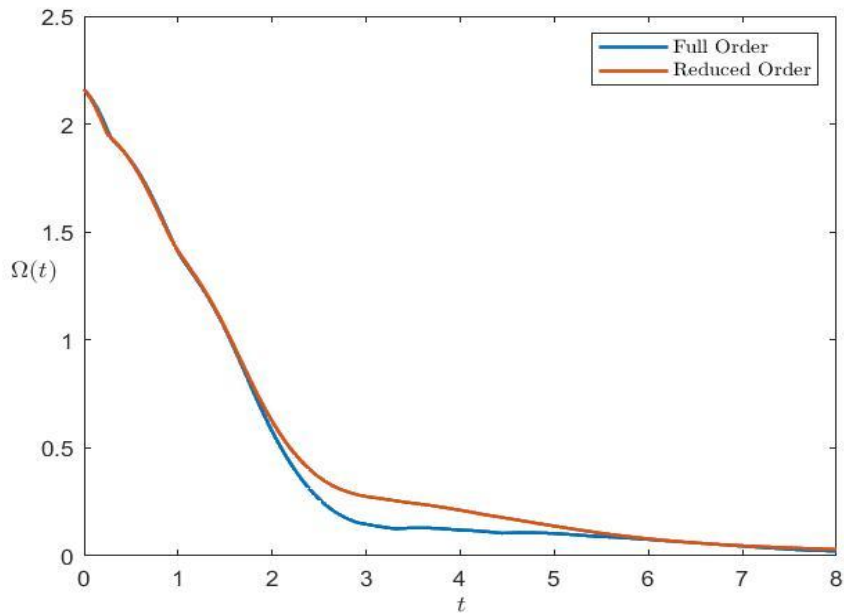

**Fig. 2:** Evolution of the function $\Omega$ in (6.13) for the case of the output feedback with the full-order observer (3.14) and the reduced-order observer (3.15).

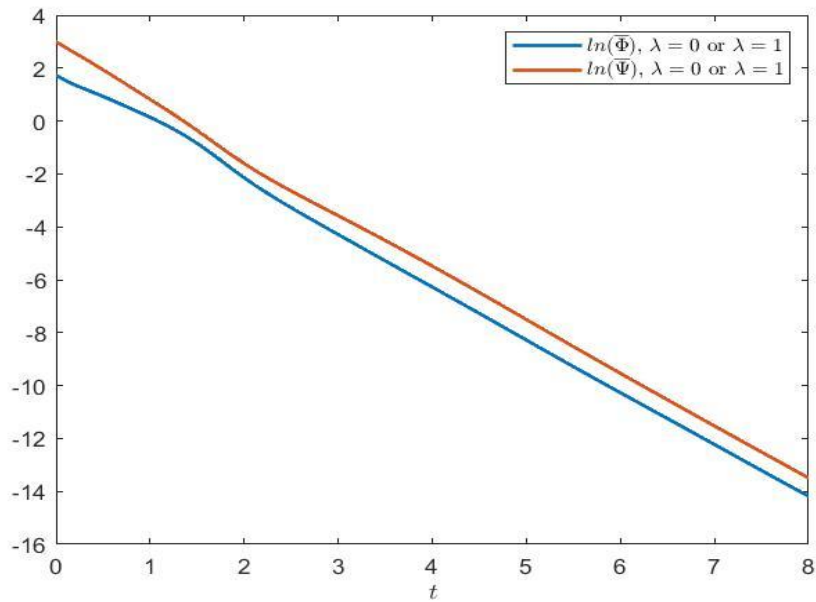

**Fig.3:** Evolution of the logarithms of the function $\bar{\Phi}$ defined by (6.11) and the function $\bar{\Psi}$ defined by (6.12).

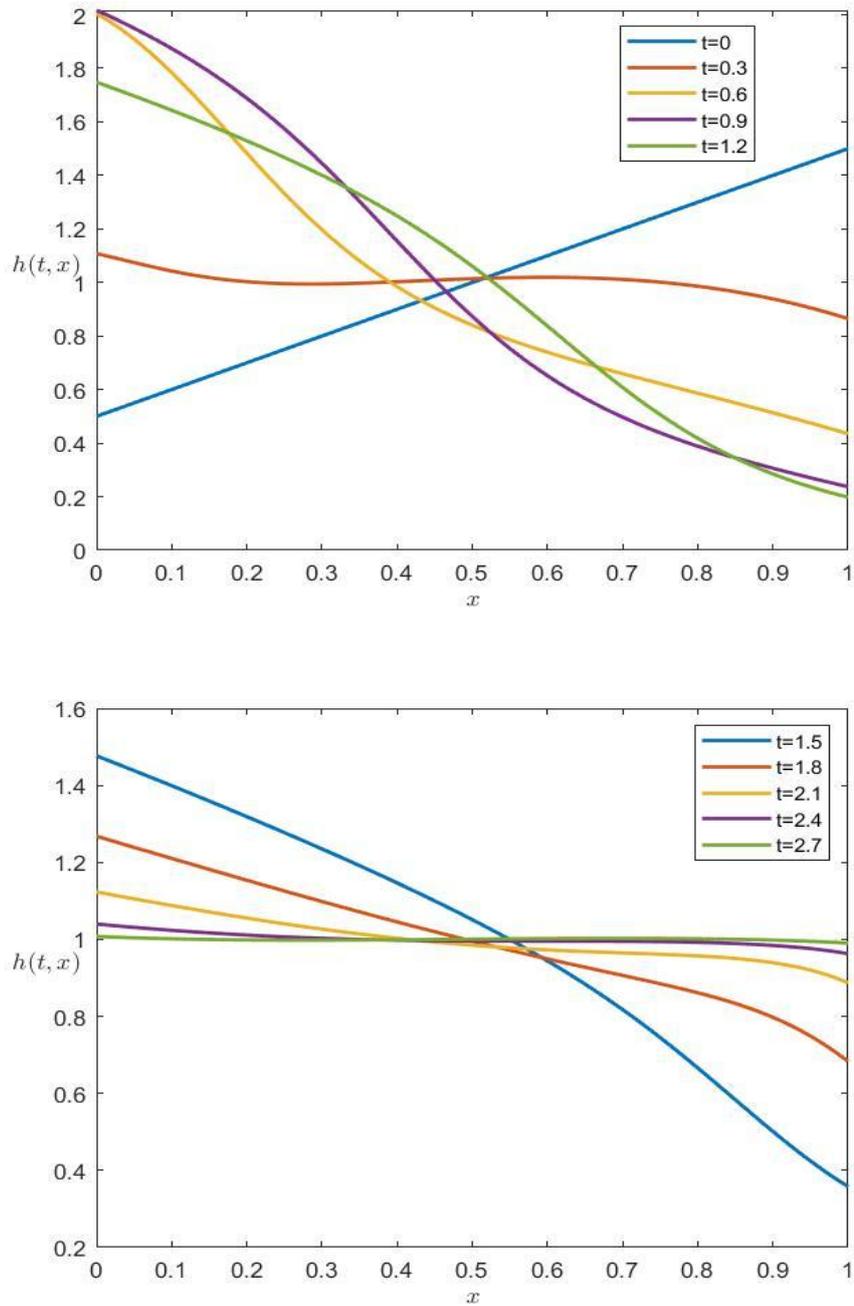

**Fig. 4:** Evolution of the liquid level for the case of the output feedback with the reduced-order observer (3.15)..

$$\dot{\Phi} = -\mu g \int_0^L h_x^2(x)dx - \mu \int_0^L h(x)v_x^2(x)dx - qk^3\xi^2$$
$$+ qk(w+k\xi)^2 - \sigma\left(2\bar{z} - \mu(h(L)-h(0)) - q(w+k\xi)\right)^2$$
$$-\sigma q^2\left(\hat{\xi}-\xi\right)^2 - 2\sigma q^2\left(\hat{w}-w-\gamma\left(\hat{\xi}-\xi\right)\right)^2$$
$$-\beta\lambda(z-\bar{z})^2 - \beta\lambda(z-\bar{z})\int_0^L h(x)v(x)dx \qquad (7.3)$$
$$-2\lambda\sigma(z-\bar{z})\left(2\bar{z} - \mu(h(L)-h(0)) - q(w+k\xi)\right)$$
$$+\sigma q\left(\hat{w}-w-\gamma\left(\hat{\xi}-\xi\right)\right)\left(2\bar{z} - \mu(h(L)-h(0)) - q(w+k\xi)\right)$$
$$+2\sigma(1-\lambda)\left(2\bar{z} - \mu(h(L)-h(0)) - q(w+k\xi)\right)\int_0^L h(x)v(x)dx$$

Using the Young inequalities

$$-(z-\bar{z})\int_0^L h(x)v(x)dx \leq \frac{1}{2}(z-\bar{z})^2 + \frac{1}{2}\left(\int_0^L h(x)v(x)dx\right)^2 \qquad (7.4)$$

$$-(z-\bar{z})\left(2\bar{z} - \mu(h(L)-h(0)) - q(w+k\xi)\right)$$
$$\leq (z-\bar{z})^2 + \frac{1}{4}\left(2\bar{z} - \mu(h(L)-h(0)) - q(w+k\xi)\right)^2 \qquad (7.5)$$

$$\left(\hat{w}-w-\gamma\left(\hat{\xi}-\xi\right)\right)\left(2\bar{z} - \mu(h(L)-h(0)) - q(w+k\xi)\right)$$
$$\leq q\left(\hat{w}-w-\gamma\left(\hat{\xi}-\xi\right)\right)^2 + \frac{1}{4q}\left(2\bar{z} - \mu(h(L)-h(0)) - q(w+k\xi)\right)^2 \qquad (7.6)$$

$$\left(2\bar{z} - \mu(h(L)-h(0)) - q(w+k\xi)\right)\int_0^L h(x)v(x)dx$$
$$\leq \frac{1}{4}\left(2\bar{z} - \mu(h(L)-h(0)) - q(w+k\xi)\right)^2 + \left(\int_0^L h(x)v(x)dx\right)^2 \qquad (7.7)$$

we obtain from (7.3) the following differential inequality that holds for all $t>0$:

$$\dot{\Phi} \leq -\mu g \int_0^L h_x^2(x)dx - \mu \int_0^L h(x)v_x^2(x)dx - qk^3\xi^2$$
$$+ qk(w+k\xi)^2 - \frac{\sigma}{4}\left(2\bar{z} - \mu(h(L)-h(0)) - q(w+k\xi)\right)^2$$
$$-\sigma q^2\left(\hat{\xi}-\xi\right)^2 - \sigma q^2\left(\hat{w}-w-\gamma\left(\hat{\xi}-\xi\right)\right)^2 \qquad (7.8)$$
$$-\lambda\left(\frac{\beta}{2} - 2\sigma\right)(z-\bar{z})^2 + \left(\frac{\beta\lambda}{2} + 2\sigma(1-\lambda)\right)\left(\int_0^L h(x)v(x)dx\right)^2$$

Using the inequality

$$2\left(2\int_0^L h(x)v(x)dx + \mu(h(L)-h(0))\right)q(w+k\xi)$$
$$\leq \frac{q^2}{2}(w+k\xi)^2 + 2\left(2\int_0^L h(x)v(x)dx + \mu(h(L)-h(0))\right)^2 \qquad (7.9)$$

we obtain from (7.8) for all $t>0$:

$$\dot{\Phi} \leq -\mu g \int_0^L h_x^2(x)dx - \mu \int_0^L h(x)v_x^2(x)dx - qk^3\xi^2$$
$$+\frac{\sigma}{4}\left(2\int_0^L h(x)v(x)dx + \mu(h(L)-h(0))\right)^2$$
$$-q\left(\frac{\sigma q}{8}-k\right)(w+k\xi)^2 - \sigma q^2\left(\hat{\xi}-\xi\right)^2 - \sigma q^2\left(\hat{w}-w-\gamma\left(\hat{\xi}-\xi\right)\right)^2 \qquad (7.10)$$
$$-\lambda\left(\frac{\beta}{2}-2\sigma\right)(z-\bar{z})^2 + \left(\frac{\beta\lambda}{2}+2\sigma(1-\lambda)\right)\left(\int_0^L h(x)v(x)dx\right)^2$$

Using the inequality

$$2\mu(h(L)-h(0))\left(2\int_0^L h(x)v(x)dx\right) \leq 4\left(\int_0^L h(x)v(x)dx\right)^2 + \mu^2(h(L)-h(0))^2 \qquad (7.11)$$

and the facts that $\beta > 4\sigma$, $\lambda \in [0,1]$, $h(L)-h(0) = \int_0^L h_x(x)dx$, we obtain from (7.10) for all $t>0$:

$$\dot{\Phi} \leq -\mu g \int_0^L h_x^2(x)dx - \mu \int_0^L h(x)v_x^2(x)dx - qk^3\xi^2 + \frac{\sigma\mu^2}{2}\left(\int_0^L h_x(x)dx\right)^2$$
$$-q\left(\frac{\sigma q}{8}-k\right)(w+k\xi)^2 - \sigma q^2\left(\hat{\xi}-\xi\right)^2 - \sigma q^2\left(\hat{w}-w-\gamma\left(\hat{\xi}-\xi\right)\right)^2 \qquad (7.12)$$
$$-\lambda\left(\frac{\beta}{2}-2\sigma\right)(z-\bar{z})^2 + \left(\frac{\beta}{2}+2\sigma\right)\left(\int_0^L h(x)v(x)dx\right)^2$$

The Cauchy-Schwarz inequality in conjunction with (2.12) implies the following inequalities

$$\left(\int_0^L h_x(x)dx\right)^2 \leq L\int_0^L h_x^2(x)dx$$
$$\left(\int_0^L h(x)v(x)dx\right)^2 \leq m\int_0^L h(x)v^2(x)dx \qquad (7.13)$$

which combined with (7.12) give for all $t>0$:

$$\dot{\Phi} \leq -\mu\left(g-\frac{\sigma\mu L}{2}\right)\int_0^L h_x^2(x)dx - \mu\int_0^L h(x)v_x^2(x)dx - qk^3\xi^2$$
$$-q\left(\frac{\sigma q}{8}-k\right)(w+k\xi)^2 - \sigma q^2\left(\hat{\xi}-\xi\right)^2 - \sigma q^2\left(\hat{w}-w-\gamma\left(\hat{\xi}-\xi\right)\right)^2 \qquad (7.14)$$
$$-\lambda\left(\frac{\beta}{2}-2\sigma\right)(z-\bar{z})^2 + \left(\frac{\beta}{2}+2\sigma\right)m\int_0^L h(x)v^2(x)dx$$

Define for each $t \geq 0$:
$$h_{\max}(t) := \max_{0 \leq x \leq L}(h(t,x)), \quad h_{\min}(t) := \min_{0 \leq x \leq L}(h(t,x)) > 0 \tag{7.15}$$

The fact that $h_{\min}(t) > 0$ follows from the fact that $h \in C^1([0,+\infty) \times [0,L]; (0,+\infty))$ (which implies that for every $t \geq 0$, $h[t]$ is a continuous positive function on $[0,L]$). Applying Wirtinger's inequality and using definition (7.15), we obtain for each $t \geq 0$:

$$\int_0^L h(x)v^2(x)dx \leq h_{\max} \int_0^L v^2(x)dx \leq \frac{L^2 h_{\max}}{\pi^2} \int_0^L v_x^2(x)dx \leq \frac{L^2 h_{\max}}{\pi^2 h_{\min}} \int_0^L h(x) v_x^2(x)dx \tag{7.16}$$

Combining (7.14) with (7.16), we get for all $t > 0$:

$$\dot{\Phi} \leq -\mu\left(g - \frac{\sigma \mu L}{2}\right)\int_0^L h_x^2(x)dx - (\beta + 4\sigma)\frac{mL^2}{2\pi^2}\left(\frac{2\pi^2 \mu}{(\beta+4\sigma)mL^2} - \frac{h_{\max}}{h_{\min}}\right)\int_0^L h(x)v_x^2(x)dx$$
$$-qk^3\xi^2 - q\left(\frac{\sigma q}{8} - k\right)(w+k\xi)^2 - \sigma q^2(\hat{\xi}-\xi)^2 - \sigma q^2\left(\hat{w}-w-\gamma(\hat{\xi}-\xi)\right)^2 \tag{7.17}$$
$$-\lambda\left(\frac{\beta}{2} - 2\sigma\right)(z - \bar{z})^2$$

Inequality (7.17) in conjunction with the facts that $\beta > 4\sigma$, $\frac{\sigma q}{8} - k > 0$, $g - \frac{\sigma \mu L}{2} > 0$ (recall (4.1)) shows that the following implication holds:

"If $t > 0$ and $\frac{2\mu\pi^2}{(\beta+4\sigma)mL^2} > \frac{h_{\max}(t)}{h_{\min}(t)}$ then $\dot{\Phi} \leq 0$" \tag{7.18}

The fact that $h \in C^1([0,+\infty) \times [0,L]; (0,+\infty))$ implies that both mappings $t \to h_{\max}(t)$, $t \to h_{\min}(t) > 0$, defined by (7.15) are continuous (see Proposition 2.9 on page 21 in [14]). Consequently, the mapping $t \to \frac{h_{\max}(t)}{h_{\min}(t)} \geq 1$ is continuous with $\frac{h_{\max}(0)}{h_{\min}(0)} \leq \frac{G^{-1}(cr)}{G^{-1}(-cr)}$. The latter inequality is a consequence of the fact that (recall definition (4.3))

$$V(\xi(0), w(0), h[0], v[0]) \leq \Phi(\xi(0), w(0), h[0], v[0], \hat{\xi}(0), \hat{w}(0), z(0)) \leq r$$

and (3.13). Since $\frac{2\mu\pi^2}{(\beta+4\sigma)mL^2} > \frac{G^{-1}(cr)}{G^{-1}(-cr)}$ (recall (4.1)) it follows that $\frac{h_{\max}(0)}{h_{\min}(0)} < \frac{2\mu\pi^2}{(\beta+4\sigma)mL^2}$. Therefore, by continuity there exists $T > 0$ such that $\frac{h_{\max}(t)}{h_{\min}(t)} < \frac{2\mu\pi^2}{(\beta+4\sigma)mL^2}$ for all $t \in [0,T)$.

We next prove by contradiction that $\frac{h_{\max}(t)}{h_{\min}(t)} < \frac{2\mu\pi^2}{(\beta+4\sigma)mL^2}$ for all $t \geq 0$. Assume the contrary, i.e. that there exists $t \geq 0$ such that $\frac{h_{\max}(t)}{h_{\min}(t)} \geq \frac{2\mu\pi^2}{(\beta+4\sigma)mL^2}$. Therefore, the set $A := \left\{ t \geq 0 : \frac{h_{\max}(t)}{h_{\min}(t)} \geq \frac{2\mu\pi^2}{(\beta+4\sigma)mL^2} \right\}$ is non-empty and bounded from below. Thus we can define $t^* = \inf(A)$. By definition, it holds that $t^* \geq T > 0$ and $\frac{h_{\max}(t)}{h_{\min}(t)} < \frac{2\mu\pi^2}{(\beta + 4\sigma)mL^2}$ for all $t \in [0, t^*)$. By continuity of the mapping $t \to \frac{h_{\max}(t)}{h_{\min}(t)}$ we obtain that $\frac{h_{\max}(t^*)}{h_{\min}(t^*)} = \frac{2\mu\pi^2}{(\beta + 4\sigma)mL^2}$. Moreover, since $\frac{h_{\max}(t)}{h_{\min}(t)} < \frac{2\mu\pi^2}{(\beta + 4\sigma)mL^2}$ for all $t \in [0, t^*)$, it follows from implication (7.18) that $\dot{\Phi} \leq 0$ for all $t \in (0, t^*)$. By continuity of the mapping $t \to \Phi\left(\xi(t), w(t), h[t], v[t], \hat{\xi}(t), \hat{w}(t), z(t)\right)$, we obtain that

$$\Phi\left(\xi(t^*), w(t^*), h[t^*], v[t^*], \hat{\xi}(t^*), \hat{w}(t^*), z(t^*)\right)$$
$$\leq \Phi\left(\xi(0), w(0), h[0], v[0], \hat{\xi}(0), \hat{w}(0), z(0)\right) \leq r$$

On the other hand, the above inequality and (4.3), (3.13) imply that $\frac{h_{\max}(t^*)}{h_{\min}(t^*)} \leq \frac{G^{-1}(cr)}{G^{-1}(-cr)} < \frac{2\mu\pi^2}{(\beta + 4\sigma)mL^2}$ which contradicts the equation $\frac{h_{\max}(t^*)}{h_{\min}(t^*)} = \frac{2\mu\pi^2}{(\beta + 4\sigma)mL^2}$.

Since $\frac{h_{\max}(t)}{h_{\min}(t)} < \frac{2\mu\pi^2}{(\beta + 4\sigma)mL^2}$ for all $t \geq 0$, we conclude from implication (7.18) that $\dot{\Phi} \leq 0$ for all $t > 0$. By continuity of the mapping $t \to \Phi\left(\xi(t), w(t), h[t], v[t], \hat{\xi}(t), \hat{w}(t), z(t)\right)$ and definition (4.3), we obtain that

$$V\left(\xi(t), w(t), h[t], v[t]\right)$$
$$\leq \Phi\left(\xi(t), w(t), h[t], v[t], \hat{\xi}(t), \hat{w}(t), z(t)\right) \qquad (7.19)$$
$$\leq \Phi\left(\xi(0), w(0), h[0], v[0], \hat{\xi}(0), \hat{w}(0), z(0)\right) \leq r < R$$

for all $t \geq 0$. Hence, $(\xi(t), w(t), h[t], v[t]) \in X$ for all $t \geq 0$ (recall definitions (3.9), (3.11)). Moreover, (3.13) implies that $\frac{h_{\max}(t)}{h_{\min}(t)} \leq \frac{G^{-1}(cr)}{G^{-1}(-cr)}$ for all $t \geq 0$. The previous inequality in conjunction with (7.17) gives for all $t > 0$:

$$\dot{\Phi} \leq -\omega\left(\int_0^L h_x^2(x)dx + \int_0^L h(x)v_x^2(x)dx + \xi^2 + \frac{\sigma q^2}{2\gamma}\left(\hat{\xi} - \xi\right)^2 + \frac{\sigma q^2}{\gamma}\left(\hat{w} - w - \gamma(\hat{\xi} - \xi)\right)^2 + (w + k\xi)^2 + \frac{\lambda}{2}(z - \bar{z})^2\right) \quad (7.20)$$

where

$$\omega := \min\left(\mu\left(g - \frac{\sigma\mu L}{2}\right), \frac{(\beta + 4\sigma)mL^2}{2\pi^2}\left(\frac{2\mu\pi^2}{(\beta + 4\sigma)mL^2} - \frac{G^{-1}(cr)}{G^{-1}(-cr)}\right), qk^3, \gamma, q\left(\frac{\sigma q}{8} - k\right), \beta - 4\sigma\right) \quad (7.21)$$

It follows from Lemma 3 in [21], definition (4.3) and (7.20) that the following estimate holds for all $t > 0$:

$$\frac{d}{dt}\Phi\left(\xi(t), w(t), h[t], v[t], \hat{\xi}(t), \hat{w}(t), z(t)\right) \leq -\frac{\omega \Phi\left(\xi(t), w(t), h[t], v[t], \hat{\xi}(t), \hat{w}(t), z(t)\right)}{\max\left(1, \Gamma\left(V\left(\xi(t), w(t), h[t], v[t]\right)\right)\right)} \quad (7.22)$$

where $\Gamma:[0,R) \to (0,+\infty)$ is a non-decreasing function. Since $\Gamma:[0,R) \to (0,+\infty)$ is non-decreasing and since $V(\xi(t), w(t), h[t], v[t]) \leq r$ for all $t \geq 0$ (recall (7.19)), we obtain from (7.22) the following estimate for all $t > 0$:

$$\frac{d}{dt}\Phi\left(\xi(t), w(t), h[t], v[t], \hat{\xi}(t), \hat{w}(t), z(t)\right) \leq -\frac{\omega}{\max(1, \Gamma(r))}\Phi\left(\xi(t), w(t), h[t], v[t], \hat{\xi}(t), \hat{w}(t), z(t)\right) \quad (7.23)$$

By continuity of the mapping $t \to \Phi\left(\xi(t), w(t), h[t], v[t], \hat{\xi}(t), \hat{w}(t), z(t)\right)$, the differential inequality (7.23) implies the following estimate for all $t \geq 0$:

$$\Phi\left(\xi(t), w(t), h[t], v[t], \hat{\xi}(t), \hat{w}(t), z(t)\right) \leq \exp\left(-\frac{\omega t}{\max(1, \Gamma(r))}\right)\Phi\left(\xi(0), w(0), h[0], v[0], \hat{\xi}(0), \hat{w}(0), z(0)\right) \quad (7.24)$$

Estimate (7.24) in conjunction with Lemma 4 in [21] and the inequalities

$$\frac{\sigma q^2}{\gamma \max(4, 4\gamma^2 + 1)}\left(\left(\hat{\xi}(t) - \xi(t)\right)^2 + \left(\hat{w}(t) - w(t)\right)^2\right)$$

$$\leq \frac{\sigma q^2}{2\gamma}\left(\hat{\xi}(t) - \xi(t)\right)^2 + \frac{\sigma q^2}{\gamma}\left(\hat{w}(t) - w(t) - \gamma\left(\hat{\xi}(t) - \xi(t)\right)\right)^2$$

$$\leq \sigma q^2\left(\frac{1}{\gamma} + 1 + \gamma\right)\left(\left(\hat{\xi}(t) - \xi(t)\right)^2 + \left(\hat{w}(t) - w(t)\right)^2\right)$$

implies the following estimate for all $t \geq 0$:

$$\tilde{\Phi}(t) \leq \exp\left(-\frac{\omega t}{\max(1, \Gamma(r))}\right) \frac{\max\left(G_2(V(\xi(0), w(0), h[0], v[0])), \sigma q^2\left(\frac{1}{\gamma} + 1 + \gamma\right), \frac{1}{2}\right)}{\min\left(\frac{1}{G_1(V(\xi(t), w(t), h[t], v[t]))}, \frac{\sigma q^2}{\gamma \max(4, 4\gamma^2 + 1)}, \frac{1}{2}\right)} \tilde{\Phi}(0) \quad (7.25)$$

where $G_i:[0,R) \to (0,+\infty)$, $i = 1, 2$, are non-decreasing functions and

$$\tilde{\Phi}(t) := \left\|(\xi(t), w(t), h[t] - h^*\chi_{[0,L]}, v[t])\right\|_X^2 + \left(\hat{\xi}(t) - \xi(t)\right)^2 + \left(\hat{w}(t) - w(t)\right)^2$$

$$+ \lambda\left(z(t) - \int_0^L h(t,x)v(t,x)dx - \mu(h(t,L) - h(t,0))\right)^2 \quad (7.26)$$

Since $G_i:[0,R) \to (0,+\infty)$, $i = 1, 2$, are non-decreasing functions and since $V(\xi(t), w(t), h[t], v[t]) \leq r$ for all $t \geq 0$ (recall (7.19)), we obtain from (7.25) the following estimate for all $t \geq 0$:

$$\tilde{\Phi}(t) \leq \exp\left(-\frac{\omega t}{\max(1, \Gamma(r))}\right) \frac{\max\left(G_2(r), \sigma q^2\left(\frac{1}{\gamma} + 1 + \gamma\right), \frac{1}{2}\right)}{\min\left(\frac{1}{G_1(r)}, \frac{\sigma q^2}{\gamma \max(4, 4\gamma^2 + 1)}, \frac{1}{2}\right)} \tilde{\Phi}(0) \quad (7.27)$$

Estimate (4.4) is a direct consequence of estimate (7.27) and definition (7.26). The proof is complete. ◁

**Proof of Theorem 2:** Let $r \in [0, R)$ be given (arbitrary) and let constants $\sigma, q, k, \gamma > 0$, $\lambda \in [0,1]$ for which (4.1) holds (but otherwise arbitrary) be given. Consider a classical solution of the PDE-ODE system (2.7)-(2.10), (2.12) with (3.15), (3.18), (4.5) that satisfies $\Psi(\xi(0), w(0), h[0], v[0], \zeta(0), z(0)) \leq r$.

We define
$$\bar{\zeta} := w - \gamma \xi \tag{7.28}$$

It follows from (2.7) and definition (7.28) that the following differential equation holds for all $t > 0$:
$$\frac{d\bar{\zeta}}{dt} = -\gamma \bar{\zeta} - \gamma^2 \xi - f \tag{7.29}$$

A direct consequence of (3.15) and (7.29) is the following differential equation for all $t > 0$:
$$\frac{d}{dt}\left(\zeta - \bar{\zeta}\right) = -\gamma \left(\zeta - \bar{\zeta}\right) \tag{7.30}$$

Using Lemma 2 in [21], (2.7), (3.19), (7.30) and definitions (3.16), (4.6), (7.28) we obtain for all $t > 0$:
$$\begin{aligned}\dot{\Psi} = &-\mu g \int_0^L h_x^2(x)dx - \mu \int_0^L h(x) v_x^2(x)dx - qk^3 \xi^2 \\ &+ qk(\bar{\zeta} + (\gamma + k)\xi)^2 + f\left(2\bar{z} - \mu(h(L) - h(0)) - q(\bar{\zeta} + (\gamma + k)\xi)\right) \\ &- 2\sigma q^2 \left(\zeta - \bar{\zeta}\right) - \beta \lambda (z - \bar{z})^2 - \beta \lambda (z - \bar{z}) \int_0^L h(x)v(x)dx.\end{aligned} \tag{7.31}$$

Equation (4.5) and definitions (3.16) and (7.28) imply that
$$\begin{aligned}f = &-\sigma \left(2\bar{z} - \mu(h(L) - h(0)) - q(\bar{\zeta} + (\gamma + k)\xi)\right) \\ &- 2\lambda \sigma (z - \bar{z}) + \sigma q (\zeta - \bar{\zeta}) + 2\sigma(1 - \lambda) \int_0^L h(x)v(x)dx\end{aligned} \tag{7.32}$$

for all $t > 0$. Combining equations (7.31), (7.32), we obtain for all $t > 0$:
$$\begin{aligned}\dot{\Psi} = &-\mu g \int_0^L h_x^2(x)dx - \mu \int_0^L h(x) v_x^2(x)dx - qk^3 \xi^2 \\ &+ qk(\bar{\zeta} + (\gamma + k)\xi)^2 - 2\sigma q^2 \left(\zeta - \bar{\zeta}\right)^2 - \beta \lambda (z - \bar{z})^2 - \beta \lambda (z - \bar{z}) \int_0^L h(x)v(x)dx \\ &+ \sigma q \left(\zeta - \bar{\zeta}\right)\left(2\bar{z} - \mu(h(L) - h(0)) - q(\bar{\zeta} + (\gamma + k)\xi)\right) \\ &- 2\lambda \sigma (z - \bar{z})\left(2\bar{z} - \mu(h(L) - h(0)) - q(\bar{\zeta} + (\gamma + k)\xi)\right) \\ &+ 2\sigma(1 - \lambda)\left(\int_0^L h(x)v(x)dx\right)\left(2\bar{z} - \mu(h(L) - h(0)) - q(\bar{\zeta} + (\gamma + k)\xi)\right) \\ &- \sigma \left(2\bar{z} - \mu(h(L) - h(0)) - q(\bar{\zeta} + (\gamma + k)\xi)\right)^2.\end{aligned} \tag{7.33}$$

Using inequality (7.4), inequalities (7.5), (7.7) with $w = \bar{\zeta} + \gamma \xi$ (recall (7.28)), and the following inequality

$$\left(\zeta-\bar{\zeta}\right)\left(2\bar{z}-\mu(h(L)-h(0))-q(\bar{\zeta}+(\gamma+k)\xi)\right) \leq$$
$$q\left(\zeta-\bar{\zeta}\right)^2 + \frac{1}{4q}\left(2\bar{z}-\mu(h(L)-h(0))-q(\bar{\zeta}+(\gamma+k)\xi)\right)^2 \tag{7.34}$$

we obtain from (7.33) the following differential inequality that holds for all $t>0$:

$$\begin{aligned}
\dot{\Psi} \leq &-\mu g \int_0^L h_x^2(x)dx - \mu \int_0^L h(x)v_x^2(x)dx - qk^3\xi^2 \\
&+ qk(\bar{\zeta}+(\gamma+k)\xi)^2 - \sigma q^2\left(\zeta-w+\gamma\xi\right)^2 \\
&- \lambda\left(\frac{\beta}{2}-2\sigma\right)(z-\bar{z})^2 + \left(\frac{\beta\lambda}{2}+2\sigma(1-\lambda)\right)\left(\int_0^L h(x)v(x)dx\right)^2 \\
&- \frac{\sigma}{4}\left(2\bar{z}-\mu(h(L)-h(0))-q(\bar{\zeta}+(\gamma+k)\xi)\right)^2
\end{aligned} \tag{7.35}$$

Using inequality (7.9) with $w=\bar{\zeta}+\gamma\xi$ (recall (7.28)) along with definitions (7.28), (3.16), we obtain from (7.35) for all $t>0$:

$$\begin{aligned}
\dot{\Psi} \leq &-\mu g \int_0^L h_x^2(x)dx - \mu \int_0^L h(x)v_x^2(x)dx - qk^3\xi^2 \\
&- q\left(\frac{\sigma q}{8}-k\right)(w+k\xi)^2 - \sigma q^2\left(\zeta-w+\gamma\xi\right)^2 \\
&- \lambda\left(\frac{\beta}{2}-2\sigma\right)(z-\bar{z})^2 + \left(\frac{\beta\lambda}{2}+2\sigma(1-\lambda)\right)\left(\int_0^L h(x)v(x)dx\right)^2 \\
&+ \frac{\sigma}{4}\left(2\int_0^L h(x)v(x)dx + \mu(h(L)-h(0))\right)^2
\end{aligned} \tag{7.36}$$

Using inequality (7.11) and the facts that $\beta>4\sigma$, $\lambda\in[0,1]$, $h(L)-h(0)=\int_0^L h_x(x)dx$, we obtain from (7.36) for all $t>0$:

$$\begin{aligned}
\dot{\Psi} \leq &-\mu g \int_0^L h_x^2(x)dx - \mu \int_0^L h(x)v_x^2(x)dx - qk^3\xi^2 \\
&- q\left(\frac{\sigma q}{8}-k\right)(w+k\xi)^2 - \sigma q^2\left(\zeta-w+\gamma\xi\right)^2 + \frac{\sigma\mu^2}{2}\left(\int_0^L h_x(x)dx\right)^2 \\
&- \lambda\left(\frac{\beta}{2}-2\sigma\right)(z-\bar{z})^2 + \left(\frac{\beta}{2}+2\sigma\right)\left(\int_0^L h(x)v(x)dx\right)^2
\end{aligned} \tag{7.37}$$

Inequalities (7.13) combined with (7.37) give for all $t>0$:

$$\begin{aligned}
\dot{\Psi} \leq &-\mu\left(g-\frac{\sigma\mu L}{2}\right)\int_0^L h_x^2(x)dx - \mu \int_0^L h(x)v_x^2(x)dx - qk^3\xi^2 \\
&- q\left(\frac{\sigma q}{8}-k\right)(w+k\xi)^2 - \sigma q^2\left(\zeta-w+\gamma\xi\right)^2 \\
&- \lambda\left(\frac{\beta}{2}-2\sigma\right)(z-\bar{z})^2 + \left(\frac{\beta}{2}+2\sigma\right)m\int_0^L h(x)v^2(x)dx
\end{aligned} \tag{7.38}$$

Using definition (7.15) and the inequalities (7.16) and (7.38), we get for all $t > 0$:

$$\dot{\Psi} \leq -\mu\left(g - \frac{\sigma\mu L}{2}\right)\int_0^L h_x^2(x)dx - \frac{(\beta+4\sigma)mL^2}{2\pi^2}\left(\frac{2\pi^2\mu}{(\beta+4\sigma)mL^2} - \frac{h_{max}}{h_{min}}\right)\int_0^L h(x)v_x^2(x)dx$$
$$-qk^3\xi^2 - q\left(\frac{\sigma q}{8} - k\right)(w+k\xi)^2 - \sigma q^2(\zeta - w + \gamma\xi)^2 - \lambda\left(\frac{\beta}{2} - 2\sigma\right)(z-\bar{z})^2 \quad (7.39)$$

Exploiting Lemma 3 in [21], (7.39) and definition (4.6) and using similar arguments as in the proof of Theorem 1, we are in a position to get the following estimates for all $t \geq 0$:

$$V(\xi(t), w(t), h[t], v[t]) \leq \Psi(\xi(t), w(t), h[t], v[t], \zeta(t), z(t)) \leq \Psi(\xi(0), w(0), h[0], v[0], \zeta(0), z(0)) \leq r < R \quad (7.40)$$

$$\Psi(\xi(t), w(t), h[t], v[t], \zeta(t), z(t)) \leq \exp\left(-\frac{\omega t}{\max(1, \Gamma(r))}\right)\Psi(\xi(0), w(0), h[0], v[0], \zeta(0), z(0)) \quad (7.41)$$

where $\omega$ is given by (7.21) and $\Gamma : [0, R) \to (0, +\infty)$ is a non-decreasing function. Estimate (7.41) in conjunction with Lemma 4 in [21] and definition (4.6) implies the following estimate for all $t \geq 0$:

$$\tilde{\Psi}(t) \leq \exp\left(-\frac{\omega t}{\max(1, \Gamma(r))}\right) \frac{\max\left(G_2(V(\xi(0), w(0), h[0], v[0])), \frac{\sigma q^2}{\gamma}, \frac{1}{2}\right)}{\min\left(\frac{1}{G_1(V(\xi(t), w(t), h[t], v[t]))}, \frac{\sigma q^2}{\gamma}, \frac{1}{2}\right)} \tilde{\Psi}(0) \quad (7.42)$$

where $G_i : [0, R) \to (0, +\infty)$, $i = 1, 2$, are non-decreasing functions and

$$\tilde{\Psi}(t) := \left\|(\xi(t), w(t), h[t] - h^*\chi_{[0,L]}, v[t])\right\|_X^2 + (\zeta(t) - w(t) + \gamma\xi(t))^2$$
$$+ \lambda\left(z(t) - \int_0^L h(t,x)v(t,x)dx - \mu(h(t,L) - h(t,0))\right)^2 \quad (7.43)$$

Since $G_i : [0, R) \to (0, +\infty)$, $i = 1, 2$, are non-decreasing functions and since $V(\xi(t), w(t), h[t], v[t]) \leq r$ for all $t \geq 0$ (recall (7.40)), we obtain from (7.42) the following estimate for all $t \geq 0$:

$$\tilde{\Psi}(t) \leq \exp\left(-\frac{\omega t}{\max(1, \Gamma(r))}\right) \frac{\max\left(G_2(r), \frac{\sigma q^2}{\gamma}, \frac{1}{2}\right)}{\min\left(\frac{1}{G_1(r)}, \frac{\sigma q^2}{\gamma}, \frac{1}{2}\right)} \tilde{\Psi}(0) \quad (7.44)$$

Estimate (4.7) is a direct consequence of estimate (7.44) and definition (7.43). The proof is complete. ◁

**Proof of Theorem 3:** Let constants $\sigma, q, k, \gamma, \beta > 0$, $\lambda \in [0,1]$ be given. Let $T > 0$ be a given constant and let $(\xi(t), w(t), h[t], v[t], \hat{\xi}(t), \hat{w}(t), z(t)) \in X \times \Re^3$ with $t \in [0, T]$ be a classical smooth solution of the PDE-ODE system (2.7)-(2.10), (2.12), (3.14), (3.18), (4.2) with $L = g = h^* = 1$. Let $n \geq 2$ be an integer and define $\delta x = 1/n$. Let $\delta t \in (0,T)$ be a given constant and consider $h_i^+, v_i^+$ for $i = 0,\ldots,n$, $w^+, \xi^+, z^+, \hat{\xi}^+, \hat{w}^+$ given by the difference scheme (6.1)-(6.6), (6.9) with $h_i = h(0, i\delta x)$, $v_i = v(0, i\delta x)$, for $i = 0,\ldots,n$, $w = w(0)$, $\xi = \xi(0)$, $\hat{\xi} = \hat{\xi}(0)$, $\hat{w} = \hat{w}(0)$ and $z = z(0)$.

Define:
$$u(t,x) = \ln(h(t,x)), \text{ for } t \in [0,T], \, x \in [0,1] \tag{7.45}$$

$$u_i = u(0, i\delta x), \text{ for } i = 0,\ldots,n \tag{7.46}$$

Using (6.1), (7.46), the following relations for $i = 1,\ldots,n-1$

$$u_x(0, i\delta x) = \frac{u_{i+1} - u_{i-1}}{2\delta x} - \frac{1}{2\delta x} \int_{i\delta x}^{(i+1)\delta x} \int_{i\delta x}^{s} \int_{i\delta x}^{l} u_{xxx}(0,\eta) d\eta dl ds - \frac{1}{2\delta x} \int_{(i-1)\delta x}^{i\delta x} \int_{s}^{i\delta x} \int_{l}^{i\delta x} u_{xxx}(0,\eta) d\eta dl ds \tag{7.47}$$

$$v_x(0, i\delta x) = \frac{v_{i+1} - v_{i-1}}{2\delta x} - \frac{1}{2\delta x} \int_{i\delta x}^{(i+1)\delta x} \int_{i\delta x}^{s} \int_{i\delta x}^{l} v_{xxx}(0,\eta) d\eta dl ds - \frac{1}{2\delta x} \int_{(i-1)\delta x}^{i\delta x} \int_{s}^{i\delta x} \int_{l}^{i\delta x} v_{xxx}(0,\eta) d\eta dl ds \tag{7.48}$$

and the fact that $u_t(t,x) = -v(t,x)u_x(t,x) - v_x(t,x)$ for $t \in [0,T]$, $x \in [0,1]$ (a consequence of (2.8) and definition (7.45)), we get for $i = 1,\ldots,n-1$

$$\begin{aligned}
u(\delta t, i\delta x) &= u_i + \int_0^{\delta t} u_t(s, i\delta x) ds = u_i - \int_0^{\delta t} v(s, i\delta x) u_x(s, i\delta x) ds - \int_0^{\delta t} v_x(s, i\delta x) ds \\
&= u_i - \delta t v_x(0, i\delta x) - \delta t v_i u_x(0, i\delta x) \\
&\quad - \int_0^{\delta t} \left( v(s, i\delta x) u_x(s, i\delta x) - v_i u_x(0, i\delta x) \right) ds - \int_0^{\delta t} \left( v_x(s, i\delta x) - v_x(0, i\delta x) \right) ds \\
&= \ln(h_i^+) + \frac{\delta t}{2\delta x} \left( \int_{i\delta x}^{(i+1)\delta x} \int_{i\delta x}^{s} \int_{i\delta x}^{l} v_{xxx}(0,\eta) d\eta dl ds + \int_{(i-1)\delta x}^{i\delta x} \int_{s}^{i\delta x} \int_{l}^{i\delta x} v_{xxx}(0,\eta) d\eta dl ds \right) \\
&\quad + \frac{\delta t}{2\delta x} v_i \left( \int_{i\delta x}^{(i+1)\delta x} \int_{i\delta x}^{s} \int_{i\delta x}^{l} u_{xxx}(0,\eta) d\eta dl ds + \int_{(i-1)\delta x}^{i\delta x} \int_{s}^{i\delta x} \int_{l}^{i\delta x} u_{xxx}(0,\eta) d\eta dl ds \right) \\
&\quad - \int_0^{\delta t} \int_0^{s} \left( v_{xt}(l, i\delta x) + v_t(l, i\delta x) u_x(l, i\delta x) + v(l, i\delta x) u_{xt}(l, i\delta x) \right) dl ds
\end{aligned} \tag{7.49}$$

Consequently, we get from (7.45) and (7.49) for all $i = 1,\ldots,n-1$:

$$\left| \ln(h(\delta t, i\delta x)) - \ln(h_i^+) \right| \leq S(\delta t) \left( (\delta x)^2 + (\delta t) \right) \tag{7.50}$$

where
$$S := \|v[0]\|_\infty \|u_{xxx}[0]\|_\infty + \|v_{xxx}[0]\|_\infty + \frac{1}{2} \max_{0 \leq s \leq T} \left( \|v_{xt}[s]\|_\infty \right) + \frac{1}{2} \max_{0 \leq s \leq T} \left( \|v_t[s]u_x[s] + v[s]u_{xt}[s]\|_\infty \right) \tag{7.51}$$

Using (6.1), (7.46), the fact that $\delta x = 1/n$ (which in conjunction with (2.10) imply that $v_0 = v_n = 0$), the following relations

$$v_x(0,0) = \frac{4v_1 - v_2}{2\delta x} - \frac{2}{\delta x}\int_0^{\delta x}\int_0^s\int_0^l v_{xxx}(0,\eta)d\eta dl ds + \frac{1}{2\delta x}\int_0^{2\delta x}\int_0^s\int_0^l v_{xxx}(0,\eta)d\eta dl ds$$

$$v_x(0,1) = \frac{v_{n-2} - 4v_{n-1}}{2\delta x} - \frac{2}{\delta x}\int_{1-\delta x}^1\int_s^1\int_l^1 v_{xxx}(0,\eta)d\eta dl ds + \frac{1}{2\delta x}\int_{1-2\delta x}^1\int_s^1\int_l^1 v_{xxx}(0,\eta)d\eta dl ds$$

and the facts that $u_t(t,0) = -v_x(t,0)$, $u_t(t,1) = -v_x(t,1)$ for $t \in [0,T]$ (consequence of (2.8), (2.10) and definition (7.45)), we get

$$u(\delta t, 0) = u_0 + \int_0^{\delta t} u_t(s,0)ds = u_0 - \int_0^{\delta t} v_x(s,0)ds$$

$$= u_0 - \delta t v_x(0,0) - \int_0^{\delta t}(v_x(s,0) - v_x(0,0))ds \qquad (7.52)$$

$$= \ln(h_0^+) - \frac{\delta t}{2\delta x}\left(\int_0^{2\delta x}\int_0^s\int_0^l v_{xxx}(0,\eta)d\eta dl ds - 4\int_0^{\delta x}\int_0^s\int_0^l v_{xxx}(0,\eta)d\eta dl ds\right) - \int_0^{\delta t}\int_0^s v_{xt}(l,0)dl ds$$

$$u(\delta t, 1) = u_n + \int_0^{\delta t} u_t(s,1)ds = u_n - \int_0^{\delta t} v_x(s,1)ds$$

$$= u_n - \delta t v_x(0,1) - \int_0^{\delta t}(v_x(s,1) - v_x(0,1))ds \qquad (7.53)$$

$$= \ln(h_n^+) + \frac{\delta t}{2\delta x}\left(4\int_{1-\delta x}^1\int_s^1\int_l^1 v_{xxx}(0,\eta)d\eta dl ds - \int_{1-2\delta x}^1\int_s^1\int_l^1 v_{xxx}(0,\eta)d\eta dl ds\right) - \int_0^{\delta t}\int_0^s v_{xt}(l,1)dl ds$$

Consequently, we get from (7.45), (7.50), (7.51), (7.52), (7.53) and the fact that $\delta x = 1/n$:

$$\left|\ln(h(\delta t, i\delta x)) - \ln(h_i^+)\right| \leq S(\delta t)\left((\delta x)^2 + (\delta t)\right), \text{ for } i = 0,1,\ldots,n \qquad (7.54)$$

Using the inequality $|\exp(x) - \exp(y)| \leq \exp(x)\exp(|x-y|)|x-y|$ that holds for all $x, y \in \Re$ and the fact that $\left|\ln(h(\delta t, i\delta x)) - \ln(h_i^+)\right| \leq ST(1+T)$ (a consequence of (7.54) and the facts that $\delta t \in (0,T)$ and $\delta x = 1/n$), we get from (7.54):

$$\max_{i=0,\ldots,n}\left(\left|h_i^+ - h(\delta t, i\delta x)\right|\right) \leq \max_{0 \leq s \leq T}\left(\|h[s]\|_\infty\right)\exp(ST(1+T))S(\delta t)\left((\delta x)^2 + (\delta t)\right) \qquad (7.55)$$

Using (3.14) and (2.7), we obtain the following formulas:

$$\hat{\xi}(\delta t) - \xi(\delta t) =$$
$$\exp(-\gamma \delta t)\left[\left(\hat{\xi}(0) - \xi(0)\right)\cos\left(\frac{\delta t}{\sqrt{2}}\right) + \sqrt{2}\left((\hat{w}(0) - w(0)) - \gamma\left(\hat{\xi}(0) - \xi(0)\right)\right)\sin\left(\frac{\delta t}{\sqrt{2}}\right)\right] \qquad (7.56)$$

$$\hat{w}(\delta t) - w(\delta t) =$$
$$\exp(-\gamma \delta t)\left[(\hat{w}(0) - w(0))\cos\left(\frac{\delta t}{\sqrt{2}}\right) + \gamma\sqrt{2}\left(\hat{w}(0) - w(0) - \left(\frac{1}{2\gamma} + \gamma\right)\left(\hat{\xi}(0) - \xi(0)\right)\right)\sin\left(\frac{\delta t}{\sqrt{2}}\right)\right] \qquad (7.57)$$

Consequently, we get from (6.5) and (7.56), (7.57):

$$\left|\hat{\xi}^+ - \hat{\xi}(\delta t)\right| = \left|\xi^+ - \xi(\delta t)\right|, \quad \left|\hat{w}^+ - \hat{w}(\delta t)\right| = \left|w^+ - w(\delta t)\right| \tag{7.58}$$

Using (2.1) we obtain the following formulas:

$$\xi(\delta t) = \xi(0) + w(0)\delta t + \frac{(\delta t)^2}{2} f(0) + \int_0^{\delta t}\int_0^s\int_0^l \dot{f}(\eta)d\eta dl ds$$

$$w(\delta t) = w(0) - f(0)\delta t - \int_0^{\delta t}\int_0^s \dot{f}(l)dl ds \tag{7.59}$$

Exploiting (7.59), (6.5) and the fact that $\delta t \in (0,T)$ we get:

$$\left|\xi(\delta t) - \xi^+\right| \leq \frac{(\delta t)^3}{6} \max_{0 \leq s \leq T}\left(\left|\dot{f}(s)\right|\right) \leq T\frac{(\delta t)^2}{6} \max_{0 \leq s \leq T}\left(\left|\dot{f}(s)\right|\right)$$

$$\left|w(\delta t) - w^+\right| \leq \frac{(\delta t)^2}{2} \max_{0 \leq s \leq T}\left(\left|\dot{f}(s)\right|\right) \tag{7.60}$$

Therefore, we get from (7.58), (7.60):

$$\left|\xi^+ - \xi(\delta t)\right| + \left|w^+ - w(\delta t)\right| + \left|\hat{\xi}^+ - \hat{\xi}(\delta t)\right| + \left|\hat{w}^+ - \hat{w}(\delta t)\right|$$
$$= 2\left|\xi^+ - \xi(\delta t)\right| + 2\left|w^+ - w(\delta t)\right| \tag{7.61}$$
$$\leq (\delta t)^2 \left(\frac{T}{3} + 1\right) \max_{0 \leq s \leq T}\left(\left|\dot{f}(s)\right|\right)$$

Using (6.2), (7.45), (7.46), (7.47), (7.48), the fact that $v_{i+1} - v_{i-1} = \int_{(i-1)\delta x}^{(i+1)\delta x} v_x(0,s)ds$ for $i = 1,\ldots,n-1$, the following relations for $i = 1,\ldots,n-1$

$$h_x(0,i\delta x) = \frac{h_{i+1} - h_{i-1}}{2\delta x} - \frac{1}{2\delta x}\int_{i\delta x}^{(i+1)\delta x}\int_{i\delta x}^s\int_{i\delta x}^l h_{xxx}(0,\eta)d\eta dl ds - \frac{1}{2\delta x}\int_{(i-1)\delta x}^{i\delta x}\int_s^{i\delta x}\int_l^{i\delta x} h_{xxx}(0,\eta)d\eta dl ds \tag{7.62}$$

$$v_{xx}(0,i\delta x) = \frac{v_{i+1} - 2v_i + v_{i-1}}{(\delta x)^2} - \frac{1}{(\delta x)^2} P_i \tag{7.63}$$

where

$$P_i := \int_{i\delta x}^{(i+1)\delta x}\int_{i\delta x}^s\int_{i\delta x}^l\int_{i\delta x}^\eta v_{xxxx}(0,\bar{r})d\bar{r} d\eta dl ds - \int_{(i-1)\delta x}^{i\delta x}\int_s^{i\delta x}\int_l^{i\delta x}\int_\eta^{i\delta x} v_{xxxx}(0,\bar{r})d\bar{r} d\eta dl ds \tag{7.64}$$

we get from (2.8), (2.9) and (7.45) for all $i = 1,\ldots,n-1$:

$$v(\delta t, i\delta x) = v_i + \delta t\, v_t(0, i\delta x) + \int_0^{\delta t} \left( v_t(s, i\delta x) - v_t(0, i\delta x) \right) ds$$

$$= v_i + \delta t \left( \mu v_{xx}(0, i\delta x) + \mu u_x(0, i\delta x) v_x(0, i\delta x) + f - v_i v_x(0, i\delta x) - h_x(0, i\delta x) \right) + \int_0^{\delta t}\int_0^s v_{tt}(l, i\delta x)\,dl\,ds$$

$$= v_i^+ - \frac{\delta t}{2\delta x} v_i \left( \int_{i\delta x}^{(i+1)\delta x}\int_{i\delta x}^s\int_{i\delta x}^l v_{xxx}(0,\eta)\,d\eta\,dl\,ds + \int_{(i-1)\delta x}^{i\delta x}\int_s^{i\delta x}\int_l^{i\delta x} v_{xxx}(0,\eta)\,d\eta\,dl\,ds \right)$$

$$- \frac{\delta t\, \mu}{4(\delta x)^2} \left( \int_{(i-1)\delta x}^{(i+1)\delta x} v_x(0,s)\,ds \right) \left( \int_{i\delta x}^{(i+1)\delta x}\int_{i\delta x}^s\int_{i\delta x}^l u_{xxx}(0,\eta)\,d\eta\,dl\,ds + \int_{(i-1)\delta x}^{i\delta x}\int_s^{i\delta x}\int_l^{i\delta x} u_{xxx}(0,\eta)\,d\eta\,dl\,ds \right)$$

$$- \frac{\delta t\, \mu}{2\delta x} u_x(0, i\delta x) \left( \int_{i\delta x}^{(i+1)\delta x}\int_{i\delta x}^s\int_{i\delta x}^l v_{xxx}(0,\eta)\,d\eta\,dl\,ds + \int_{(i-1)\delta x}^{i\delta x}\int_s^{i\delta x}\int_l^{i\delta x} v_{xxx}(0,\eta)\,d\eta\,dl\,ds \right)$$

$$+ \frac{\delta t}{2\delta x} \left( \int_{i\delta x}^{(i+1)\delta x}\int_{i\delta x}^s\int_{i\delta x}^l h_{xxx}(0,\eta)\,d\eta\,dl\,ds + \int_{(i-1)\delta x}^{i\delta x}\int_s^{i\delta x}\int_l^{i\delta x} h_{xxx}(0,\eta)\,d\eta\,dl\,ds \right) - \frac{\delta t\, \mu}{(\delta x)^2} P_i + \int_0^{\delta t}\int_0^s v_{tt}(l, i\delta x)\,dl\,ds$$

(7.65)

It follows from (7.64), (7.65) that the following estimate holds for $i = 1, \ldots, n-1$:

$$\left| v(\delta t, i\delta x) - v_i^+ \right| \leq \bar{S}\, \delta t \left( (\delta x)^2 + \delta t \right) \tag{7.66}$$

where

$$\bar{S} := \frac{1}{2} \left( \max_{0 \leq s \leq T} \left( \|v_{tt}[s]\|_\infty \right) + \|v[0]\|_\infty \|v_{xxx}[0]\|_\infty + \|h_{xxx}[0]\|_\infty \right)$$

$$+ \frac{\mu}{2} \left( \|v_x[0]\|_\infty \|u_{xxx}[0]\|_\infty + \|u_x[0]\|_\infty \|v_{xxx}[0]\|_\infty + \|v_{xxxx}[0]\|_\infty \right) \tag{7.67}$$

By virtue of (2.10), (6.3) and (7.67) we conclude that

$$\max_{i=0,\ldots,n} \left( \left| v_i^+ - v(\delta t, i\delta x) \right| \right) \leq \bar{S}\, \delta t \left( \delta t + (\delta x)^2 \right) \tag{7.68}$$

Using (3.18) (with $m = L = g = 1$), the fact that $\delta x = 1/n$ and (6.8) we get:

$$z(\delta t) - z^+ = \int_0^{\delta t}\int_0^s \ddot{z}(l)\,dl\,ds \tag{7.69}$$

Consequently, (7.69) implies that

$$\left| z(\delta t) - z^+ \right| \leq \frac{(\delta t)^2}{2} \max_{0 \leq s \leq T} \left( |\ddot{z}(s)| \right) \tag{7.70}$$

Combining (7.55), (7.61), (7.68) and (7.70) we obtain estimate (6.9) with $\bar{M} = \bar{S} + \frac{1}{2}\max_{0 \leq s \leq T}\left(|\ddot{z}(s)|\right) + \left(\frac{T}{3}+1\right)\max_{0 \leq s \leq T}\left(|\dot{f}(s)|\right) + \max_{0 \leq s \leq T}\left(\|h[s]\|_\infty\right)\exp(ST(1+T))S$. The proof is complete. ◁

**Proof of Theorem 4:** Let constants $\sigma, q, k, \gamma, \beta > 0$, $\lambda \in [0,1]$ be given. Let $T > 0$ be a given constant and let $(\xi(t), w(t), h[t], v[t], \zeta(t), z(t)) \in X \times \Re^2$ with $t \in [0,T]$ be a classical smooth solution of the PDE-ODE system (2.7)-(2.10), (2.12), (3.15), (3.18), (4.5) with $L = g = h^* = 1$. Let $n \geq 2$ be an integer and define $\delta x = 1/n$. Let $\delta t \in (0,T)$ be a given constant and consider $h_i^+, v_i^+$ for $i = 0, \ldots, n$, $w^+, \xi^+, z^+, \zeta^+$ given by the difference scheme (6.1)-(6.3), (6.6)-(6.8) with $h_i = h(0, i\delta x)$, $v_i = v(0, i\delta x)$, for $i = 0, \ldots, n$, $w = w(0)$, $\xi = \xi(0)$, $\zeta = \zeta(0)$, and $z = z(0)$.

Using (3.15) and (2.7), we obtain the following formula:

$$\zeta(\delta t) = w(\delta t) - \gamma \xi(\delta t) + \exp(-\gamma \delta t)(\zeta - w + \gamma \xi) \tag{7.71}$$

Consequently, we get from (6.7) and (7.71):

$$\left|\zeta(\delta t) - \zeta^+\right| = \left|w(\delta t) - w^+ + \gamma(\xi^+ - \xi(\delta t))\right| \leq \left|w(\delta t) - w^+\right| + \gamma\left|\xi^+ - \xi(\delta t)\right| \tag{7.72}$$

Using (2.1) we obtain the formulas (7.59), (7.60). Therefore, we get from (7.72), (7.60):

$$\begin{aligned}
&\left|\xi^+ - \xi(\delta t)\right| + \left|w^+ - w(\delta t)\right| + \left|\zeta^+ - \zeta(\delta t)\right| \\
&\leq (\gamma + 1)\left|\xi^+ - \xi(\delta t)\right| + 2\left|w^+ - w(\delta t)\right| \\
&\leq (\delta t)^2 \left((\gamma + 1)\frac{T}{6} + 1\right) \max_{0 \leq s \leq T}\left(\left|\dot{f}(s)\right|\right)
\end{aligned} \tag{7.73}$$

The rest of proof of Theorem 4 is exactly the same with the proof of Theorem 3. The proof is complete. ◁

## 8. Concluding Remarks

By applying the CLF methodology we have managed to achieve output feedback stabilization results for the viscous Saint-Venant liquid-tank system. As far as we know, this is the first paper in the literature that achieves output feedback stabilization of the nonlinear viscous Saint-Venant system.

The obtained results leave some open problems which will be the topic of future research:
1) The results were applied to classical solutions and it is of interest to relax this to weak solutions. It is also an open problem to show existence/uniqueness of (weak) solutions for the closed-loop system. To this purpose, ideas utilized in [41] can be employed.
2) The construction of CLFs which can allow the derivation of stability estimates in stronger spatial norms for the liquid level and velocity profiles.

Another more demanding problem that will be studied in the future is the spill-free, slosh-free and smash-free movement of a glass of water. This problem arises when we want to move the glass of water to a position which is close to a wall. In this case, we need to control the overshoot of the glass position error in order to avoid smashing the glass on the wall.

**Acknowledgements:** The authors would like to thank Dr. Dionysios Theodosis for his valuable contribution to the simulation results.